\documentclass[1p,final]{elsarticle}

\usepackage{amssymb,mathtools}
\usepackage{siunitx}

\usepackage[acronym,nopostdot,nowarn]{glossaries}
\newacronym{ai}{AI}{artificial intelligence}
\newacronym{bspm}{BSPM}{body surface potential map}
\newacronym{cdt}{CDT}{cardiac digital twin}
\newacronym{cep}{CEP}{cardiac electrophysiology}
\newacronym{csp}{CSP}{conduction system pacing}
\newacronym{ct}{CT}{X-ray computed tomography}
\newacronym[longplural={conduction velocities}]{cv}{CV}{conduction velocity}
\newacronym{cvd}{CVD}{cardiovascular disease}
\newacronym{dt}{DT}{digital twin}
\newacronym{eam}{EAM}{electro-anatomical mapping}
\newacronym{ecdt}{eCDT}{electrophysiological cardiac digital twin}
\newacronym{ecg}{ECG}{electrocardiogram}
\newacronym{ecgi}{ECGi}{electrocardiographic imaging}
\newacronym{egm}{EGM}{electrogram}
\newacronym{ep}{EP}{electrophysiology}
\newacronym{fe}{FE}{finite element}
\newacronym{gpu}{GPU}{graphical processing unit}
\newacronym{gt}{GT}{ground truth}
\newacronym{hps}{HPS}{His-Purkinje system}
\newacronym{isct}{ISCT}{\emph{in silico} clinical trial}
\newacronym{lat}{LAT}{local activation time}
\newacronym{lv}{LV}{left ventricle}
\newacronym{pmj}{PMJ}{Purkinje myocardial junction}
\newacronym{ml}{ML}{machine learning}
\newacronym{mri}{MRI}{magnetic resonance imaging}
\newacronym{qrs}{QRS}{QRS complex}
\newacronym{roi}{ROI}{region of influence}
\newacronym{rmsd}{RMSD}{root mean squared deviation}
\newacronym{rmse}{RMSE}{root mean squared error}
\newacronym{rv}{RV}{right ventricle}
\newacronym{uvc}{UVC}{universal ventricular coordinate}
\newacronym{vsc}{VSC}{Vienna scientific cluster}
\newacronym{vt}{VT}{ventricular tachycardia}
\newacronym{wct}{WCT}{Wilson's central terminal}

\newif\ifarxiv
\arxivtrue  

\ifarxiv
\newcommand{\revaat}[1]{#1}
\else
\usepackage{xcolor}
\newcommand{\revaat}[1]{\textcolor{red}{#1}}

\usepackage[mathlines]{lineno}
\linenumbers
\fi

\usepackage[hidelinks]{hyperref}
\newcommand{\abs}[1]{\left\lvert#1\right\rvert}
\newcommand{\diff}[1]{\ensuremath{\operatorname{d}\!{#1}}}

\newcommand{\vx}[0]{\boldsymbol{x}}
\newcommand{\vf}[0]{\boldsymbol{f}}
\newcommand{\vs}[0]{\boldsymbol{s}}
\newcommand{\vn}[0]{\boldsymbol{n}}

\def\geodesicBP{Geodesic-BP}

\newcommand{\Vmrest}[0]{\ensuremath{V_0}}
\newcommand{\Vmexc}[0]{\ensuremath{V_1}}
\newcommand{\V}[0]{\ensuremath{\mathcal{V}}}
\newcommand{\Ge}[0]{\ensuremath{\mathbf{G}_e}}
\newcommand{\Gi}[0]{\ensuremath{\mathbf{G}_i}}
\newcommand{\G}[0]{\ensuremath{\mathbf{G}}}
\newcommand{\TD}[0]{\ensuremath{\mathbb{T}}}
\newcommand{\Z}[1]{\ensuremath{Z_{#1}}}

\newcommand{\taumean}[0]{\ensuremath{\tau_\mu}}
\newcommand{\taustd}[0]{\ensuremath{\tau_\sigma}}
\newcommand{\Vmean}[0]{\ensuremath{\V_\mu}}
\newcommand{\Vstd}[0]{\ensuremath{\V_\sigma}}

\newcommand{\distecg}[0]{\ensuremath{\operatorname{dist}_{\V}}}

\newcommand{\distlat}[0]{\ensuremath{\operatorname{dist}_{\tau}}}
\newcommand{\distlatij}[0]{\ensuremath{\distlat (\tau_i, \tau_j)}}
\newcommand{\distpb}[0]{\ensuremath{\operatorname{dist}_{\phi}}}
\newcommand{\distpbij}[0]{\ensuremath{\distpb (\phi_i, \phi_j)}}

\title{Accurate and Efficient Cardiac Digital Twin from Surface ECGs:
       Insights into Identifiability of Ventricular Conduction System}

\author[1,2]{Thomas Grandits\corref{equal}}
\ead{thomas.grandits@uni-graz.at}

\author[3,4]{Karli Gillette\corref{equal}}
\ead{karli.gillette@medunigraz.at}

\author[3,4]{Gernot Plank\corref{cor1}}
\ead{gernot.plank@medunigraz.at}

\author[2,5]{Simone Pezzuto\corref{cor1}}
\ead{simone.pezzuto@unitn.it}

\cortext[cor1]{Corresponding authors.}
\cortext[equal]{These authors contributed equally to this work.}

\address[1]{Department of Mathematics and Scientific Computing, University of Graz, Heinrichstraße 36, 8010 Graz, Austria}
\address[2]{Euler Institute, Università della Svizzera italiana, Via Buffi 13, 6900 Lugano, Switzerland}
\address[3]{Gottfried Schatz Research Center: Medical Physics and Biophysics, Medical University of Graz, Neue Stiftingtalstraße 6, 8010 Graz, Austria}
\address[4]{BioTechMed-Graz, Mozartgasse 12/II, 8010 Graz, Austria}
\address[5]{Department of Mathematics, University of Trento, Via Sommarive 14, 38123 Trento, Italy}

\begin{document}

\begin{abstract}
    Digital twins for cardiac electrophysiology are an enabling technology for precision cardiology. Current forward models are advanced enough to simulate the cardiac electric activity under different pathophysiological conditions and accurately replicate clinical signals like torso electrocardiograms (ECGs). In this work, we address the challenge of matching subject-specific QRS complexes using anatomically accurate, physiologically grounded cardiac digital twins.
    By fitting the initial conditions of a cardiac propagation model, our non-invasive method predicts activation patterns during sinus rhythm. 
    For the first time, we demonstrate that distinct activation maps can generate identical surface ECGs. To address this non-uniqueness, we introduce a physiological prior based on the distribution of Purkinje-muscle junctions. Additionally, we develop a digital twin ensemble for probabilistic inference of cardiac activation.
    Our approach marks a significant advancement in the calibration of cardiac digital twins and enhances their credibility for clinical application.    

        
\end{abstract}

\maketitle

\section{Introduction}
\label{sec:introduction}

\Gls{cdt} technology aims to build a 1:1 \emph{in silico} replica of a subject's heart~\cite{niederer2021scaling}.
Built upon advanced biophysical simulations, \glspl{cdt} are informed with constant or periodic real-time updates 
from observations~\cite{laubenbacher2024digital},
to accurately track the physiological state of a patient's heart. 
Based on the mechanistic nature of the used modeling technologies, 
a critical element in the \gls{cdt} paradigm is the inherent assumption 
that virtual and physical \revaat{hearts} are closely intertwined
such that any stimulus or perturbation leads to the same emergent response, 
whether this occurs in the real or virtual space.
\glspl{cdt} fulfilling this assumption offer immense transformative potential for personalizing medical care,
as they provide safe and effective means of assessing a patient's cardiac function
to support diagnosis, stratification and optimal planning  of therapies ~\cite{corral_acero_digital_2020,cluitmans2024digital}. 

Realizing this promising vision of building high fidelity \glspl{cdt} at scale
poses a number of technological challenges. 
Complex modeling workflows must be implemented comprising two distinct stages
-- referred to as \emph{anatomical} and \emph{functional twinning} --
dedicated to building anatomically accurate models of the heart from medical images,
and, subsequently, calibrating a vast number of parameters of the model 
to assimilate electrophysiological function between model and patient heart, respectively.
Significant progress has been made in terms of automating anatomical twinning 
to build patient heart anatomies at scale ~\cite{crozier2016image}.
Less developed and markedly more challenging are functional twinning methodologies.
The vast majority of simulation studies using models labelled as \glspl{cdt}
refrain from functional twinning.
Rather, the same average parameters are used across the board for all patients,
and like-for-like comparisons with clinically observable signals such as \glspl{egm} or the \gls{ecg}
are avoided as these would reveal substantial discrepancies between virtual model and physical reality 
\cite{sung2022fat,bishop24:_stochastic}.
Recent efforts to develop enhanced functional twinning technologies 
show promise to lift current restrictions in the calibration of models.
A key advance was the introduction of more efficient methods for evaluating the biophysical model 
that now facilitate cardiac \gls{ep} simulations at the organ scale with real-time performance~\cite{pezzuto_evaluation_2017,neic2017efficient}.
These have been used in proof-of-concept studies to calibrate models of patient hearts 
to non-invasive clinical \gls{ecg} data using both sampling~\cite{pezzuto_reconstruction_2021,gillette2021framework,camps2024digital,alvarez_barrientos_probabilistic_2025} and optimization approaches~\cite{grandits_geasi_2021,grandits_digital_2023}.
As such, the grand challenge at the very core of functional twinning remains --
to fit the model's high dimensional parameter space from limited clinical, ideally non-invasively acquired, data
to accurately reproduce clinical recordings like-for-like with high fidelity in a robust, unique and scalable manner. 



A main objective is to infer the electrical activation sequence of the ventricles 
-- the main pumping chambers of the heart -- from a patient's surface \gls{ecg},
and to replicate the sequence in a physiologically accurate model of electrical wavefront propagation
and the associated \gls{ecg} (specifically, the QRS complex)~\cite{pezzuto_reconstruction_2021,grandits_geasi_2021,gillette2022personalized,li_deep_2022,grandits_digital_2023,camps2024digital}.
The ventricular activation sequence and its reflection in the QRS complex observed on the body surface is driven by \revaat{the} ventricular conduction system referred to as \gls{hps} that initiates the activation of the ventricles.
The \gls{hps} comprises the bundle of His 
that receives atrial signals via the atrio-ventricular node, 
a left and right bundle branch that splits into fascicles attached to the Purkinje network,
a fast-conducting network of Purkinje fibers that permeates both the endocardial and sub-endocardial tissue~\cite{spach_electrical_1963,myerburg1972physiology,vigmond_2016_modeling}. 
The Purkinje network is coupled to the ventricular myocardium at terminal junctions, referred to as \glspl{pmj},
\revaat{which} serve as the initiation sites of earliest ventricular activation and, thus, 
play a pivotal role in determining the ventricular activation pattern and the surface \gls{ecg}. 
However, the structure of the \gls{hps} as well as the distribution of \glspl{pmj} is highly patient-specific 
and, generally, inaccessible by non-invasive methods~\cite{pullan_inverse_2010}, 
and, even with invasive advanced mapping methodologies~\cite{palamara2014computational}, only observable with limited accuracy.
Owing to its importance as the pivotal structure controlling ventricular activation, and as a target for the treatment of conduction disturbances~\cite{sharma2018permanent},
a method to unveil its structure and properties 
is among the most fundamental challenges in cardiac electrophysiology.
However, inferring an extremely high dimensional object such as the ventricular activation sequence 
from very limited, sparse, and noisy data such as the \gls{ecg}, constitutes an ill-posed problem.
As such, it is still unclear 
as to which extent and accuracy the ventricular activation sequence can actually be inferred 
solely from the standard 12-lead \gls{ecg}~\cite{grandits2022smoothness}. Consequently, the possibility that different ventricular activation maps  could yield the same QRS complexes in the 12-lead \gls{ecg} must be considered~\cite{franzone_mathematical_2014,li_solving_2025}.

Fast \gls{ep} simulators, based on the Eikonal and reaction-Eikonal models, offer a way to address the identifiability challenge of inferring -- ideally uniquely -- 
the ventricular  activation sequence~\cite{franzone_spreading_1993,pezzuto_evaluation_2017,neic2017efficient}. 
These simulators allow for the rapid exploration of the parameter space 
to identify candidate \gls{pmj} sets by comparing simulated against recorded \glspl{ecg} ~\cite{gillette2021framework,camps2024digital}. 
Sampling-based methods have successfully identified ventricular activation patterns 
from ECG data under conditions such as focal activity~\cite{pezzuto_learning_2022,meisenzahl2024boatmap}, 
bundle branch block~\cite{pezzuto_reconstruction_2021,alvarez_barrientos_probabilistic_2025}, 
and during normal sinus rhythm~\cite{gillette2021framework,camps2024digital,cardone2016human,barber2021estimation}. 
However, given the enormous number of samples required, 
these methods are feasible only for a small number of parameters 
or when paired with fast emulators~\cite{pezzuto_learning_2022,salvador_digital_2023,camps2024digital}. 
Here, the \gls{hps} is often approximated by a small number of \glspl{pmj} combined with a fast endocardial layer~\cite{pezzuto_reconstruction_2021}.
In general though, sampling-based approaches and emulators do not scale well as the number of parameters increases~\cite{larson2019derivative},
thus they are limited in their ability to infer the ventricular conduction system in a more general setting.

\begin{figure}[!t]
       \centering
        \includegraphics[width=\linewidth]{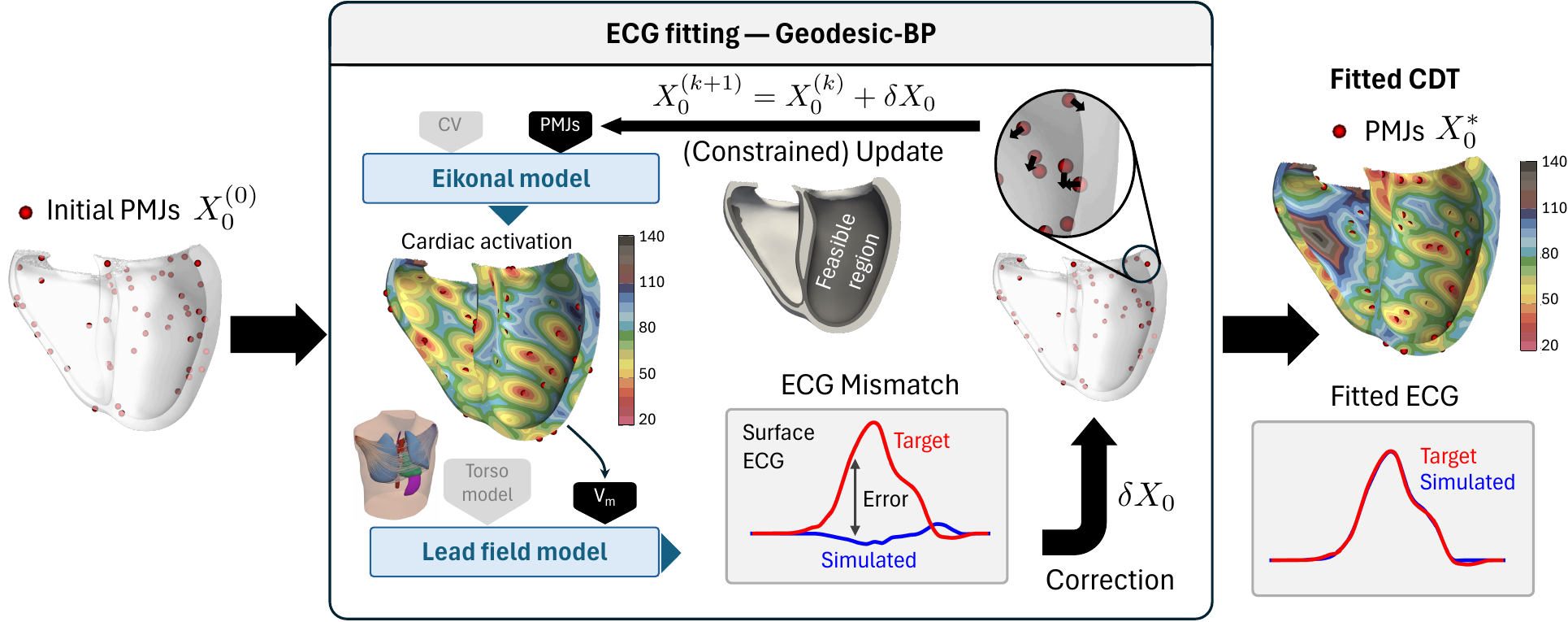}
       \caption{
       \textbf{Geodesic-BP}: fast definition of a \gls{cdt} from the surface \gls{ecg}. Geodesic-BP~\cite{grandits_digital_2023} optimizes the distribution of \glspl{pmj} until the mismatch between recorded and simulated \gls{ecg} is minimized. Physiological constraints for the \glspl{pmj} are automatically imposed during the optimization process.
       }
       \label{fig:methods_overview}
\end{figure}

In this work, we leverage \geodesicBP, a fast, gradient-based method for identifying \glspl{pmj} from \gls{ecg} data~\cite{grandits_digital_2023}. (See Fig.~\ref{fig:methods_overview} for an overview of \geodesicBP.) 
Geodesic-BP can achieve a highly accurate fit to the \gls{ecg} in under 30 minutes on a single high-end \gls{gpu}. 
\revaat{Note that while the method can also be computed on a CPU, the runtime increases to more than 6 hours.}
Being gradient-based, it scales efficiently with the number of parameters 
and guarantees convergence to a locally optimal solution. 
This makes it an ideal platform for investigating the identifiability problem. 
In general, Geodesic-BP converges toward a local minimum of the loss function, 
which measures the mismatch between recorded and simulated \gls{ecg}. 
The presence of multiple minima indicates a lack of identifiability, 
while the distribution of these minima in the parameter space offers insights into potential regularization strategies.

\revaat{We tackle this issue by fitting the \gls{ecg} of a single case several hundred times} with unprecedented accuracy, 
using different initial parameter conditions and varying the number of ECG leads, 
including up to high-density \glspl{bspm}. 
We then quantify the uncertainty in the identification process. 
Our results are quantitatively validated against high-fidelity \gls{gt} data generated \emph{in silico} 
using a different \gls{ep} simulator \cite{gillette2022personalized}. 
Finally, we propose physiologically grounded constraints for the \glspl{pmj} 
that significantly improve the identifiability of the ventricular conduction system from the 12-lead \gls{ecg}. Exact constraints from the \gls{gt} model were not utilized to elucidate how physiological assumptions of the \gls{hps} may \revaat{may play a role} in the ability to obtain an exact match in the \gls{ecg}. 

This is the first study demonstrating the feasibility of inferring the \gls{hps} at unprecedented accuracy
from a standard clinical 12-lead \gls{ecg}.
Our approach produces \revaat{highly accurate} \gls{ecg} matches with a narrow spread in uncertainty 
in the identified ventricular activation times that is well below the uncertainty of measurements
with the most advanced invasive mapping technologies.
We believe that our novel robust and scalable optimization and forward \gls{ep} approach 
is the key technology for creating credible \glspl{cdt} with demonstrable fidelity,
and, thus, for translating \glspl{cdt} to clinical applications.

\section{Methods}


\subsection{Electrophysiology model}
\label{sec:epmodel}

The forward model for the \gls{ecg} is based on the eikonal equation with the lead field method~\cite{neic2017efficient,pezzuto_evaluation_2017,gillette2021framework}.
Considering a domain $\Omega\subset\mathbb{R}^3$ representing the active ventricular myocardium, the eikonal equation provides the cardiac \gls{lat} \revaat{$\tau(\vx):\Omega\to\mathbb{R}$} given a \revaat{symmetric, positive-definite} conduction velocity tensor \revaat{$\mathbf{M}(\vx):\Omega\to\mathbb{R}^{3\times 3}$} and a set of boundary conditions $X_0$:
\begin{equation}
\left\{ \begin{aligned}
        &\sqrt{\left(\mathbf{M}(\vx) \nabla \tau(\vx) \right) \cdot \nabla \tau(\vx)} = 1, && \text{in $\Omega$}, \\
        &\tau(\vx_i) = t_i, && \text{for $(\vx_i, t_i) \in X_0$},
\end{aligned} \right.
\label{eq:anisotropic_eikonal_eq}
\end{equation}
where
\begin{align}
\mathbf{M}(\vx) &= v_f^2 \vf(\vx) \otimes \vf(\vx) 
+ v_s^2 \vs(\vx) \otimes \vs(\vx) 
+ v_n^2 \vn(\vx) \otimes \vn(\vx), \\
X_0 &= \bigl\{ (\vx_i,t_i) \colon \vx_i \in \bar\Omega, t_i \in \mathbb{R}, i = 1,\ldots,N \bigr\},
\end{align}
and $v_{\{f,s,n\}}(\vx)$ are the \glspl{cv} in the fiber~$\vf$, sheet~$\vs$ and normal~$\vn$ direction, respectively. The set $X_0$ contains the \glspl{pmj} $\vx_i$ and their activation time $t_i$, which should satisfy the following compatibility condition:
\begin{equation}
\revaat{t_i \le t_j + \delta(\vx_i,\vx_j)}, \quad \text{for all $(\vx_i,t_i),(\vx_j,t_j)\in X_0$},
\label{eq:compat}
\end{equation}
where $\delta(\vx,\boldsymbol{y})$ is the travel time between $\vx$ and $\boldsymbol{y}$ according to Eq.~\eqref{eq:anisotropic_eikonal_eq}.

The numerical approximation, previously described~\cite{grandits_digital_2023}, consists in a piecewise linear approximation of $\tau(\vx)$ on a simplicial mesh of $\Omega$. The conduction velocity tensor is assumed piecewise constant in each element of the grid. The solution is based on a fixed-point iteration scheme where $\tau(\vx)$ at each node is updated according to the neighbor values of $\tau$:
\begin{equation}
\tau^{(k+1)}(\vx_i) = \min_{\boldsymbol{y}\in\partial\omega(\vx_i)} \Bigl\{ \tau^{(k)}(\vx_i), \tau^{(k)}(\boldsymbol{y}) + \delta(\vx_i,\boldsymbol{y}) \Bigr\},
\label{eq:hopf_lax_update}
\end{equation}
where $\vx_i$ is the position of the $i$-th mesh node and $\omega(\vx_i)$ is the set of elements sharing~$\vx_i$. 
Eq.~\eqref{eq:hopf_lax_update} is applied in parallel to all nodes until convergence. The geodesic distance is computed by solving a local optimization problem on $\partial\omega(\vx_i)$ with a fixed number of FISTA algorithm~\cite{beck_fast_2009}  iterations. 
When a \gls{pmj} $(\vx_j,t_j)\in X_0$ is not a mesh node, we solve a local eikonal problem in the mesh element containing $\vx_j$, and then assign the boundary condition $\tau(\vx_k) = t_j + \delta(\vx_j,\vx_k)$ for all neighbour nodes $\vx_k$. 
Finally, the compatibility condition \eqref{eq:compat} is automatically fulfilled by the algorithm, in the sense that when violated, Eq.~\eqref{eq:anisotropic_eikonal_eq} overwrites the corresponding \gls{pmj} with the correct timing. 
In this way, the final number of active \glspl{pmj} could be lower than~$N$.
Overall, our vectorized approach is highly efficient for \glspl{gpu} and it is suitable for automatic differentation (see Sec.~\ref{sec:geodesic_bp}).


\subsection{ECG computation}
\label{sec:ecg}

The transmembrane voltage $V_\text{m}(\vx,t)$ during depolarization is reconstructed from the \gls{lat} $\tau(\vx)$ as follows:
\begin{equation}
    V_\text{m}(\vx,t; X_0) = U\bigl(t-\tau(\vx;X_0)\bigr)  
    = \Vmrest + \frac{\Vmexc - \Vmrest}{2} \left[ \tanh \left( 2\frac{t - \tau(\vx;X_0)}{\varepsilon} \right) + 1 \right],
    \label{eq:tanh_waveform}
\end{equation}
where $\Vmrest$ and $\Vmexc$ refer to resting and plateau potential. 
\revaat{This simplified ionic model can be derived by assuming the cell state to be modeled by a bistable equation~\cite[Chapter 6.2]{keener_mathematical_1998}.}
Note that we emphasize the parametric dependency of $V_\text{m}(\vx,t; X_0)$ on $X_0$.
In all experiments we chose $\Vmrest = \SI{-85}{\milli\volt}$, $\Vmexc = \SI{30}{\milli\volt}$, and $\varepsilon = \SI{1}{\milli\second}$ 
to closely approximate the upstroke of a human ventricular action potential.  

The \gls{ecg} is obtained with the lead field method~\cite{pezzuto_evaluation_2017,potse2018scalable}. 
For a single lead $l_e=1,\ldots,L$, the corresponding \gls{egm} $\V_{l_e}(t)$ is given by
\begin{equation}
    \V_{l_e}(t; X_0) = \int_\Omega \Gi(\vx) \nabla \Z{l_e}(\vx) \cdot \nabla V_\text{m}(\vx,t;X_0) \diff{\vx},
    \label{eq:integral_potential_ecg}
\end{equation}
where $Z_{l_e}(\vx)$ is the lead field corresponding to the lead $l_e$ 
and $\Gi$ is the intracellular conductivity tensor. (See Sec.~\ref{sec:lead_field} for the computation of $Z_{l_e}$.) 
\revaat{We refer to the combination of all \glspl{egm} as the \gls{ecg}.}

For a fixed time-step $\Delta t$, we define a temporal grid $t_0, t_0 + \Delta t, \ldots, t_0 + N_t \Delta t = t_\text{end}$ and compute $\V_{l_e}(t;X_0)$ as follows:
\begin{equation}
    \V_{l_e}(t_k;X_0) = \mathbf{B}_{l_e} \cdot \mathbf{V}_\text{m}(t_k;X_0),
\end{equation}
where $\mathbf{V}_\text{m}(t_k;X_0)$ is the (interpolated) transmembrane potential at time $t_k$ and $\mathbf{B}_{l_e}$ is the vector:
\begin{equation}
[\mathbf{B}_{l_e}]_j = \int_{\Omega} \Gi(\vx) \nabla \Z{l_e}(\vx) \cdot \nabla \psi_j(\vx) \diff{\vx},
\end{equation}
with $\psi_j(\vx)$ the piecewise linear hat function corresponding to the mesh node $\vx_j$, that is $\psi_j(\vx_j) = 1$ and $\psi_j(\vx_i) = 0$ for all $i\neq j$. Note that the vector $\mathbf{B}_{l_e}$ is constant and can be precomputed for each lead. 
For the temporal discretization parameter, we select $\Delta t = \SI{0.5}{\ms}$,
and $t_\text{end}$ is chosen such that $t_\text{end} > \Delta \text{QRS}$ holds,
with $\Delta \text{QRS}$ being the duration of the entire QRS complex in the 12-lead ECG.


\subsection{Model calibration via Geodesic-BP}
\label{sec:geodesic_bp}
To achieve efficient cardiac electrophysiological twinning through the surface ECG, 
we employ \revaat{Geodesic-BP}, as previously presented in~\cite{grandits_digital_2023},
and briefly summarized here.

Given the recorded \gls{ecg} $\hat{\V}_{l_e}(t)$, $l_e=1,\ldots,L$,
we can pose the problem of finding an optimal set of initiation sites $X_0$ such that the simulated \gls{ecg} $\V_{l_e}(t;X_0)$ matches in a least-squares sense the recorded one:
\begin{equation}
\min_{X_0 \in \mathcal{S}} 
\frac{1}{L}\frac{1}{|\TD|}
\sum_{l_e=1}^L \sum_{k=0}^{N_t} \Bigl( \V_{l_e}(t_k;X_0) - \hat{\V}_{l_e}(t_k) \Bigr)^2 
= \min_{X_0 \in \mathcal{S}} \mathcal{L}(X_0),
    \label{eq:inverse_ecg}
\end{equation}
\revaat{for} $\abs{\TD}=t_\text{end}-t_0$ and $\mathcal{S}$ \revaat{being} the feasible set of \glspl{pmj} and $\mathcal{L}(X_0)$ is the loss function. The feasible set allows for stricter geometrical constraints to improve identifiability and \revaat{is defined} as follows:
\begin{equation}
\mathcal{S} = S \times \{ t \ge 0 \},
\end{equation}
that is, positions $\vx_i$ of \glspl{pmj} must lie within the subdomain $S\subset\bar\Omega$, and their activation times $t_i$ must be non-negative. 
No upper bound is \revaat{imposed on} $t_i$; consequently, \glspl{pmj} with large $t_i$ are automatically deactivated by the algorithm due to the compatibility condition~\eqref{eq:compat}. The subdomain $S$ can be either: 
\begin{itemize}
\item (unrestricted case) the entire active myocardium, or
\item (restricted case) a subendocardial band, as described in Sec.~\ref{sec:constraints}.
\end{itemize}

Numerically, we minimize the loss function in Eq.~\eqref{eq:inverse_ecg} using the gradient-based optimizer ADAM~\cite{kingma_adam_2017} for $400$ iterations with a learning rate of $0.75$. Starting from an initial guess $X_0^{(0)}\subset\mathcal{S}$ for the \glspl{pmj}, the optimization algorithm generates a sequence $X_0^{(k)}, X_0^{(k+1)}, \ldots$ of \glspl{pmj} that reduces the mismatch with the targeted \gls{ecg} until convergence. We denote by $X_0^*$ the optimal set of \glspl{pmj}. The update $X_0^{(k+1)}$ is based on the gradient of the loss function, computed via backpropagation -- for this problem, backpropagation returns discrete geodesic paths from the \glspl{pmj}, which explains the name of the algorithm -- and is then projected onto $\mathcal{S}$ by taking the closest point in $S$.

\subsection{Sampling solutions}
\label{sec:sampling}

In general, Geodesic-BP will converge towards different optimal \glspl{pmj}, $X_0^*$, depending on the choice of $X^{(0)}_0$, that is the initial guess for the \glspl{pmj} for the optimization algorithm.
We empirically quantify the variability of optimized solutions as follows.
We generate 20 initial sets of \glspl{pmj} by random sampling, with each set comprising 300 \glspl{pmj} locations and timings. 
The random sampling of $X^{(0)}_0$ is as follows: each \gls{pmj} has initial position, $\vx_i$, uniformly distributed \revaat{on the boundary of the feasible domain} $\partial S$ using \texttt{trimesh}~\cite{dawson_haggerty_et_al_trimesh_2019}, and $t_i$ uniformly distributed between $t_0$ and $t_\text{end}$. Note that all sites and timings are drawn independently. By default, the number of \glspl{pmj} is set to $N=300$, unless otherwise stated.
Each of them is optimized by Geodesic-BP to minimize the mismatch to the recorded \gls{gt} \gls{ecg}.


\subsection{Geometrical constraints for PMJs}
\label{sec:constraints}

In the restricted case, \glspl{pmj} must lie inside the subendocardial band $S$, which has been chosen based on the known structure of the \gls{hps} from histo-anatomical studies. 
%
%
\begin{figure}[t]
    \centering
    \includegraphics[width=.85\linewidth, trim={0.2cm, 2.5cm 0cm 2.5cm}, clip]{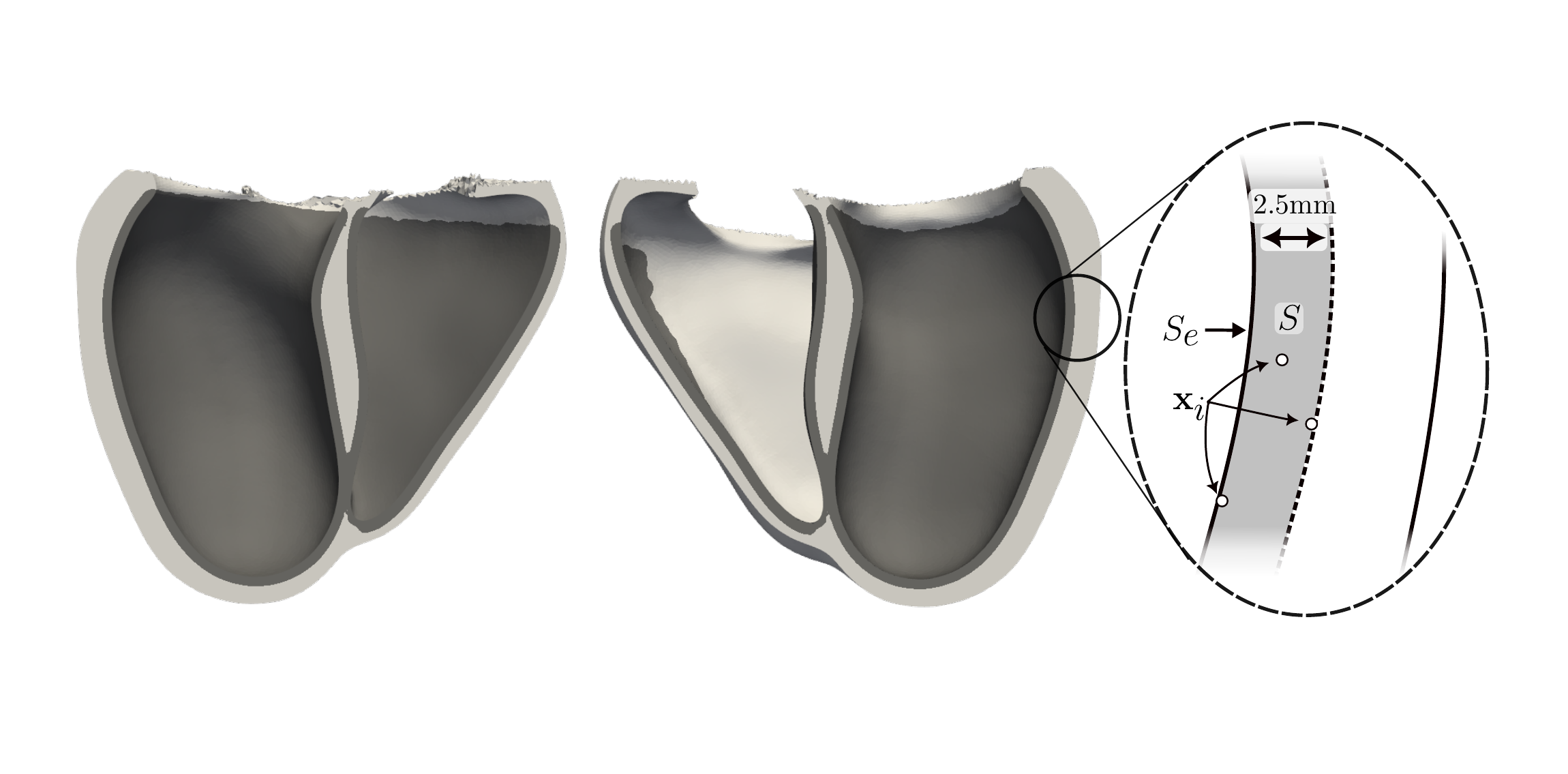}
    \caption{Endocardial surface restrictions limits the domain of possible locations for \glspl{pmj}. 
    The \glspl{pmj} are restricted to lie inside the subdomain $S$ which is spanned by a band, equidistant to the endocardium $S_e$.  
    }
    \label{fig:insilico_restrictions}
\end{figure}
%
%
%
%
As shown in Fig.~\ref{fig:insilico_restrictions}, $S$ is the region comprising all points within a given distance $d_\text{PMJ} = \SI{2.5}{\mm}$ from an endocardial surface $S_e\subset\partial\Omega$:
\begin{equation}
S = \bigl\{ \vx \in \bar{\Omega}: \operatorname{dist}(\vx,S_e)\le d_\text{PMJ} \bigr\},
\end{equation}
where $\operatorname{dist}(\vx,S_e)$ is the minimal Euclidean distance of $\vx$ from $S_e$.

To determine the extent of $S_e$, we utilized information on the coverage and density of the \gls{hps} across the endocardium. 
The \gls{hps} begins in the atrioventricular node and connects to a His-Bundle within the ventricular septum. To achieve synchronous activation of the \gls{lv} and the \gls{rv}, the \gls{hps} further splits after the His-bundle into the left and right bundle branches that separate further into 3 and 2 fascicles, respectively \cite{stephenson_2017_high, durrer_1970_total, myerburg1972physiology}. Each fascicle spreads across the subendocardium \cite{myerburg1972physiology} and has a fractal-like structure as there is a merging of Purkinje fibers in addition to bifurcations \cite{vigmond_2016_modeling}. The HPS is generally less dense in the RV than LV with limited coverage towards the base and in the inferior free wall of the RV \cite{stephenson_2017_high,myerburg1972physiology}, leading to delayed activation \cite{durrer_1970_total}. 
Therefore, the surface $S_e$: 1) lies within $\approx \SI{90}{\percent}$ of the apico-basal extent, i.e.\ excluding a rim of \SI{10}{\percent} apico-basal extent underneath the ventricular base, and 2) must not be located within the \gls{rv} inferior wall. This \revaat{is consistent} with the known restrictions from the \gls{gt} model.
In terms of intramural depth of \glspl{pmj}, histological evidence suggests 
that Purkinje fibers in the \gls{hps} reside mostly sub-endocardially, 
with a potential transmural penetration of $<30\%$ of the wall width~\cite{vigmond2016modeling}. 
To avoid compromising the generalizability of the method, the exact geometric constraints 
according to the \gls{pmj} depth underlying the \gls{gt} model (detailed later in Sec.~\ref{sec:gt_model}) were not exploited. Instead, a general depth of $\SI{2.5}{\milli \meter}$ globally over the entire endocardial domain was selected 
and deemed as permissible based on ventricular wall thicknesses in the general population \cite{sjogren1971left,prakash1978determination}.
\revaat{To note, the topology and the conduction velocity of the \gls{hps} are not available to the inverse problem in any form.}

\subsection{Pseudo-bidomain computations}
\label{sec:pseudo_bidomain}
\revaat{In order to assess the fidelity of our algorithm on the whole torso, we utilize the pseudo-bidomain equation to compute the extracellular potentials on the torso surface.}
The extracellular potentials $\phi_0(\vx,t)$ in the torso $\Omega_0\subset\mathbb{R}^3$ are based on the pseudo-bidomain equation~\cite{bishop_bidomain_2011}. 
Given $V_\text{m}(\vx,t)$ in the active myocardium, for each time $t\in[t_0,t_\text{end}]$ we solve the following elliptic problem:
\begin{equation}
\left\{ \begin{aligned}
-\operatorname{div} \bigl(\G \nabla \phi \bigr) &= \operatorname{div} \bigl(\Gi \nabla V_\text{m} \bigr), && \vx \in \Omega, \\
-\operatorname{div} \bigl(\G_0 \nabla \phi_0 \bigr) &= 0, && \vx \in \Omega_0, \\
\phi &= \phi_0, && \vx \in \Gamma_H, \\
\G \nabla \phi\cdot\vn - \G_0 \nabla \phi_0\cdot\vn &= \Gi \nabla V_\text{m}\cdot\vn, && \vx \in \Gamma_H, \\
\G_0 \nabla \phi_0\cdot\vn &= 0, && \vx \in \Gamma_T,
\end{aligned}\right.
\label{eq:pseudo_bidomain}
\end{equation}
where $\G = \Ge + \Gi$ is the bulk conductivity, \revaat{for $\Gi, \Ge$ being the intra-/extracellular conductivity tensors respectively}, $\phi(\vx,t)$ is the extracellular potential in the myocardium $\Omega$, $\Gamma_H = \bar\Omega \cap \bar\Omega_0$ is the heart-torso interface, and $\Gamma_T = \partial\Omega_0 \setminus \Gamma_H$ is the torso surface. Since the solution is defined up to a constant, we impose the following constraint~\cite{franzone_mathematical_2014}:
\begin{equation}
   \int_{\Gamma_T} \phi(\vx, t) \diff{\vx} = 0.
   \label{eq:torso_zero}
\end{equation}

The intracellular and extracellular conductivities ($\Gi$ and $\Ge$ respectively) within the ventricles were based on \citet{roberts1982effect} at values of $g_{i,f}=\SI{0.34}{\siemens \per \meter}$, $g_{e,f}=\SI{0.12}{\siemens \per \meter}$, $g_{i,s}=g_{i,n}=\SI{0.06}{\siemens \per \meter}$, and $g_{e,s}=g_{e,n}=\SI{0.08}{\siemens \per \meter}$. \revaat{The atrial conductivities were tuned to values of $g_{i,f}=\SI{1.02}{\siemens \per \meter}$, $g_{e,f}=\SI{3.72}{\siemens \per \meter}$, $g_{i,s}=g_{i,n}=\SI{0.22}{\siemens \per \meter}$, and $g_{e,s}=g_{e,n}=\SI{2.77}{\siemens \per \meter}$ according to conductivity ratios in \cite{roberts1982effect} as reported in \cite{nagel2022comparison,gillette2022personalized}.} All remaining conductivities in the torso, $\G_0$, were assigned the nominal values reported in \citet{keller2010ranking}, such that the torso, blood pools, and lungs \revaat{have} isotropic conductivities of $0.22$, $0.7$, and $\SI{0.0389}{\siemens \per \meter}$ respectively. 


Eq.~\eqref{eq:pseudo_bidomain} is numerically discretized with linear finite elements on a tetrahedral mesh of the torso (see Sec.~\ref{sec:gt_model}). The corresponding linear system is solved with the conjugate gradient method with the orthogonalization of the residual with the linear constraint~\eqref{eq:torso_zero} at the discrete level~\cite{bochev2005}.

\subsection{Lead fields}
\label{sec:lead_field}

Since the \gls{ecg} or the \gls{bspm} are sparse evaluations of the torso potential, $\phi_0(\vx,t)$, it is convenient to adopt the lead field formulation~\cite{potse2018scalable}.
The \revaat{\gls{egm}} for a given lead $l_e$, $\V_{l_e}(t)$, is generally defined as follows:
\begin{equation}
\V_{l_e}(t) = \sum_e w_e \phi_0(\vx_e,t),
\end{equation}
where $\vx_e\in \Gamma_T$ are \revaat{the electrode positions} and $w_e$ are the corresponding weights. 
The lead field $Z_{l_e}(\vx)$ is defined as the solution of the following problem:
\begin{equation}
\left\{ \begin{aligned}
-\operatorname{div} \bigl(\G \nabla Z_{l_e} \bigr) &= 0, && \vx \in \Omega, \\
-\operatorname{div} \bigl(\G_0 \nabla Z_0 \bigr) &= 0, && \vx \in \Omega_0, \\
Z_{l_e} &= Z_0, && \vx \in \Gamma_H, \\
\G \nabla Z_{l_e}\cdot\vn - \G_0 \nabla Z_0\cdot\vn &= 0, && \vx \in \Gamma_H, \\
\G_0 \nabla Z_0\cdot\vn &= \sum_e w_e \delta_{\vx_e}(\vx), && \vx \in \Gamma_T,
\end{aligned}\right.
\label{eq:lead_field}
\end{equation}
where $\delta_{\vx_e}$ is the Dirac measure on $\Gamma_T$ and $Z_0(\vx)$ is the lead field in the torso. Note that we are only interested in $Z_{l_e}(\vx)$, that is the lead field in the heart.

Formally, Eq.~\eqref{eq:lead_field} corresponds to Eq.~\eqref{eq:pseudo_bidomain} with a different right hand side. Thus, we impose the constraint in Eq.~\eqref{eq:torso_zero} to the solution, yielding to the compatibility condition $\sum_e w_e = 0$ on the weights. 
To keep the \gls{ecg} amplitude in clinical ranges \revaat{in the \gls{gt} model}, the computed lead fields were further scaled by a factor of $0.32$. \revaat{This scaling overcomes discrepancies in amplitudes stemming from cardiac sources within the biophysical model.} Conductivity values for lead field computation are the same as described earlier in Sec.~\ref{sec:pseudo_bidomain}.

\subsection{Ground truth model}
\label{sec:gt_model}

\begin{figure}[tb]
        \centering
        \includegraphics[width=.95\linewidth]{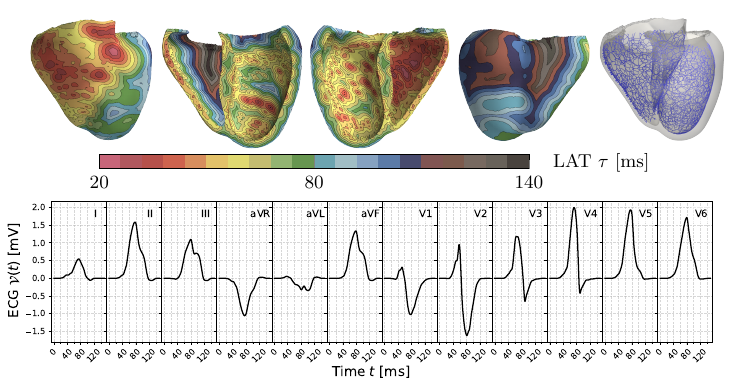}
        \caption{The \emph{in silico} \gls{gt} solution generated from the physiologically-detailed \gls{hps}.
        The \gls{lat} map $\tau$ is shown together with the \gls{hps} in blue (top).
        The bottom row shows the corresponding, calculated 12-lead \gls{ecg}.}
        \label{fig:insilico_purkinje}
\end{figure}

To test the identifiability of \gls{ecg}-based personalization, 
we utilize a state-of-the-art subject-specific biophysically detailed model
of the ventricles leveraging a topologically realistic model of the \gls{hps} 
(see Fig.~\ref{fig:insilico_purkinje}).
The model had been calibrated to replicate the ECG of this healthy subject 
with high fidelity \cite{gillette2022personalized}. 

The acquisition of the geometrical model and \gls{ecg} has been previously described elsewhere~\cite{gillette2022personalized,gillette2021framework}. 
Briefly, two magnetic resonance images were acquired sequentially in the single male subject within an approved \gls{mri} study. 
The two \gls{mri} scans were acquire sequentially: 1) a 4-stack torso \gls{mri}, 
and 2) a 3D whole-heart cardiac \gls{mri}. 
Before \gls{mri} acquisition, raw \glspl{ecg} were recorded at the ten electrode locations, $\vx_{l_e}$, in an 12-lead ECG configuration. 
Electrodes were \gls{mri}-compatible and left intact during acquisition 
to facilitate recovery in the geometric model.

A convolutional neural network was deployed for automatic cardiac segmentation 
from the 3D whole-heart \gls{mri} scan \cite{payer2018multi}, 
and the torso was segmented using semi-automatically intensity threshold approaches in \emph{Seg3D} \cite{sci_seg3D}. 
In a combined segmentation, a finite-element volumetric mesh with tetrahedral elements of approximately \SI{1253}{\micro \meter} resolution within the ventricles and up to \SI{2280}{\micro \meter} in the torso, atria, blood pools, and lungs was generated 
using \emph{meshtool}~\cite{neic2020automating}. 



Ventricular activation was facilitated through a physiologically and topologically realistic \gls{hps} (Fig.~\ref{fig:insilico_purkinje}) created assuming physiological knowledge attained experimentally within both humans and other mammalian species \cite{stephenson_2017_high, durrer_1970_total, myerburg1972physiology} more thoroughly detailed above in Sec.~\ref{sec:constraints}. 
The constructed \gls{hps} was discretized to give a spatial resolution of \SI{482}{\micro \meter}. 
There were a total of 907 \glspl{pmj} in the \gls{hps} facilitating the sites of earliest activation within the ventricles. 
\glspl{pmj} were configured to reside in up to \SI{5}{\percent} \revaat{transmurally from the endocardium for} all fascicles, apart from \revaat{up to} \SI{50}{\percent} \revaat{transmurally} within the \gls{rv} moderator band to give a realistic QRS complex~\cite{gillette2021hps}. 
To account for electrotonic loading and associated anterograde and retrograde conduction delays, 
\glspl{pmj} were coupled to the ventricular mesh using a radius of \SI{1981}{\micro \meter} 
to ensure nodal capture of 10 adjacent nodes. 
A total of 8587 nodes were selected as coupled \glspl{pmj}. 
\gls{pmj} penetration depths were computed as the closest distance to the endocardial surface. 

The \gls{cv} within the \gls{hps} was prescribed at a value of \SI{2.0}{\meter \per \second}~\cite{kassebaum1966electrophysiological}. Modification in CV to certain branches in the \gls{hps} was made to ensure earliest activation corresponded with ventricular activation as detailed in previous work~\cite{gillette2021framework}.
In contrast to the previously reported setup \cite{gillette2022personalized}, the \gls{cv} in the remaining ventricular myocardium was prescribed values of $v_f=\SI{0.61}{\meter \per \second}$ and $v_n=v_s=\SI{0.225}{\meter \per \second}$~\cite{grandits_digital_2023}. 


%
%

\subsubsection{BSPM vest generation}
\label{sec:setup_bspm}
\begin{figure}[htb]
\centering
\includegraphics[width=.75\linewidth]{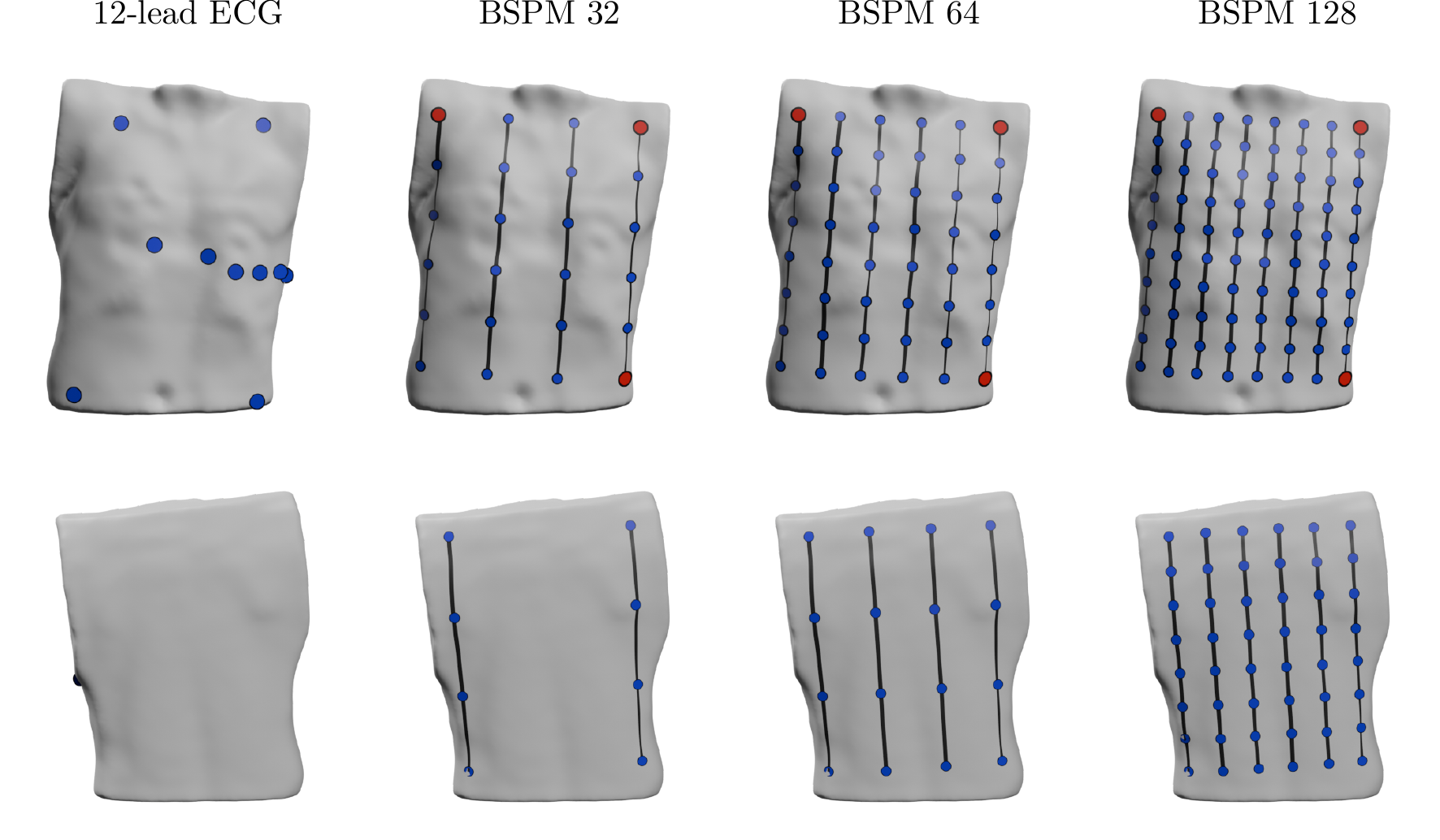}
\caption{Electrode locations of the 12-lead \gls{ecg} (left) compared to \gls{bspm} configurations comprising $32$, $64$ and $128$ electrodes from an anterior (top) and posterior (bottom) view.
For each \gls{bspm} electrode the uni-polar lead field was computed 
using the \gls{wct} as reference. Electrodes for computing the \gls{wct} are shown in red.}
\label{fig:insilico_bspms_comparison}
\end{figure}

To test the effect of denser torso electrode coverage, we equipped the anatomical torso model with an in silico \gls{bspm} vest.
These vests, commonly found in \gls{ecgi} technologies, should provide denser recordings and potentially better reconstructions~\cite{webber_technical_2023,ghanem_noninvasive_2005} and are also commonly found in explanted torso tank experiments~\cite{calder_torso_tank_2018}.
We employ a grid-based design, similar to the cited studies, with a varying amount of electrodes on the front and back of the torso.
In Fig.~\ref{fig:insilico_bspms_comparison}, we show the generated electrode positions of the artificial vest for the tested $32$, $64$ and $128$ electrodes.
For the lead field computation (Sec.~\ref{sec:lead_field}), we assume all electrodes to be unipolar w.r.t.\ the \gls{wct}, marked in red.

\subsection{Evaluation}
\label{sec:evaluation}
The ability of the Geodesic-BP method to identify the ventricular activation sequence 
from the clinical 12-lead \gls{ecg} and the goodness of fit of the associated \gls{ecg} is evaluated 
by comparing against the known \gls{gt} solution, 
and by measuring the variability of the solution as a function of initial conditions.

\subsubsection{Fidelity metrics}
The fidelity of the signal reconstruction at the body surface 
is measured for the \gls{ecg} as the root mean squared deviation between computed and \gls{gt} \glspl{ecg}
in \unit{\milli \volt} and in the interval $\TD = [t_0,t_\text{end}]$:
\begin{equation}
    \distecg (\V, \hat{\V}) = \sqrt{ \frac{1}{\abs{\TD} N_e} \sum_{l_e} \int_{\TD} \bigl(\V_{l_e}(t) - \hat{\V}_{l_e}(t) \bigr)^2 \diff{t}},
    \label{eq:ecg_dist}
\end{equation}
and, analogously, for the \gls{bspm} the difference between two \glspl{bspm} $\phi_i$,  $\phi_j$ 
over the whole torso surface $\Gamma_T$ is quantified by:
\begin{equation}
    \distpbij = \sqrt{ \frac{1}{\abs{\Gamma_T} \abs{\TD{}}} \int_{\TD} \int_{\Gamma_T} \bigl(\phi_i(\vx, t) - \phi_j(\vx, t) \bigr)^2 \diff{\vx} \diff{t}},
    \label{eq:pb_dist}
\end{equation}
where $\abs{A}$ is the measure of the set $A$.
\revaat{The difference} between two reconstructed \gls{lat} maps $\tau_i$, $\tau_j$
in \unit{\milli \second} \revaat{is} measured as:
\begin{equation}
    \distlatij = \sqrt{ \frac{1}{\abs{\Omega}} \int_\Omega \bigl(\tau_i(\vx) - \tau_j(\vx) \bigr)^2 \diff{\vx}}.
    \label{eq:sol_dist}
\end{equation}

\subsubsection{Variability and parameter relevance}
The variability between multiple solutions is quantified by computing point-wise variations of the solutions.
For this purpose, consider a set of \gls{lat} solutions $\tau_1(\vx), \ldots, \tau_N(\vx)$. 
We define the point-wise mean, $\taumean$, and standard deviation, $\taustd$, as
\begin{align}
\taumean(\vx)  &= \frac{1}{N} \sum_i \tau_i(\vx), \\
\taustd^2(\vx) &= \frac{1}{N} \sum_i \bigr( \tau_i(\vx) - \taumean(\vx) \bigl)^2.
\end{align}
The average deviation on the heart is given by $\bar{\tau}_\sigma := \frac{1}{\abs{\Omega}} \int_\Omega \taustd(\vx) \diff{\vx}$ in \unit{\milli \second}.
The mean and standard \gls{egm} for a single electrode -- $\Vmean(t)$ and $\Vstd(t)$ respectively -- are defined in the same way.

The importance of the parameters on the actual solution can also be computed numerically:
For a single \gls{pmj} $(\vx_i, t_i)$, we can compute the \gls{roi}, \revaat{that corresponds to the total influence of the \gls{pmj} on the solution $\tau$} in \unit{\cubic \milli \meter}.
\revaat{To this end, we use} the \gls{lat} derivative to compute the total influence of the \gls{pmj} on the solution $\text{ROI}_i = \int_\Omega \partial_{t_i} \tau(\vx) \diff{\vx}$.
We call \glspl{pmj} with no influence on the solution $(\text{ROI}_i = 0)$ inactive.
More information on this quantity and its derivation can be found in Sec.~\ref{app:roi}.

\subsubsection{Experimental protocols}
Different \emph{in silico} experimental scenarios were considered 
to investigate the fidelity of \gls{ecg} calibration and the identifiability of the \gls{hps} 
by varying physiological constraints as well as the density of observed data. 
The following protocols were implemented where the Geodesic-BP optimizer was employed 
to find a set of \glspl{pmj} by minimizing the mismatch to the given \gls{gt} \gls{ecg}: 
\begin{itemize}
    \item No restrictions were imposed on the location of initiation sites, 
          with permissible locations over the entire ventricular domain;
    \item Physiological constraints were imposed restricting the permissible domain 
          to subendocardial regions;
    \item Using the same physiological constraints, denser observations \revaat{beyond} the \gls{ecg} 
          were provided in the form of \glspl{bspm} of variable density, 
          using between 16 to 128 recording sites (see Fig.~\ref{fig:insilico_bspms_comparison}); 
    \item The number of initiation sites as the only hyperparameter of the optimization algorithm 
          was sampled to determine its impact upon fidelity of the reconstructed \gls{ecg}
          and the identifiability of the ventricular activation maps.
\end{itemize}

\subsection{Implementation aspects}
\label{sec:software}

Implementation aspects of the optimization algorithm Geodesic-BP 
have been reported in detail previously~\cite[II-E, Implementation aspects]{grandits_digital_2023}.
\revaat{In the present work}, the originally proposed algorithm was accelerated 
by computationally optimizing the code for solving the local problem~\cite[Algorithm 1]{grandits_digital_2023}, 
yielding a reduction in runtime of a single iteration 
from $\approx \SI{16}{\second}$ down to \SI{6}{\second}.

Geodesic-BP was implemented in Python and uses the \texttt{NumPy}~\cite{harris_array_2020} and \texttt{PyTorch}~\cite{paszke_pytorch_2019} packages. 
For the finite element computations \texttt{scikit-fem}~\cite{skfem_2020} was used. 
Simulation of the \gls{gt} electrophysiology model, \revaat{and} the verification 
of computed activation sequences and associated ECG traces were carried out 
using a biophysically detailed high fidelity reaction-eikonal model \cite{neic2017efficient} 
as provided by the established cardiac electrophysiology simulator \texttt{CARPentry}~\cite{vigmond2008solvers}.
For visualization, projection and additional 3D computations, \texttt{PyVista}, \texttt{TriMesh}, \texttt{VTK}~\cite{sullivan_pyvista_2019,dawson_haggerty_et_al_trimesh_2019,schroeder_visualization_2006} and Blender~\cite{blender_online_community_blender_2024} were also utilized.

\section{Results}
\label{sec:results}

\subsection{Identifiability in the general case}
\label{sec:uniqueness}

From a theoretical standpoint, there is no guarantee that the \gls{hps} 
governing ventricular activation can be uniquely identified from the surface \gls{ecg}. 
Optimization algorithms such as \geodesicBP\ might converge towards different local minima of the loss function, 
depending on the initial positions and timings of the \glspl{pmj}.
Here, we empirically quantify the variability of optimized solutions in the unrestricted case, 
without imposing any physiological constraints on the permissible positions of \glspl{pmj},
and by randomly sampling the initial guess for the \glspl{pmj}
in the optimization algorithm (see Section~\ref{sec:sampling}).

%
\begin{figure}[htb]
    \centering
    \includegraphics[width=\textwidth]{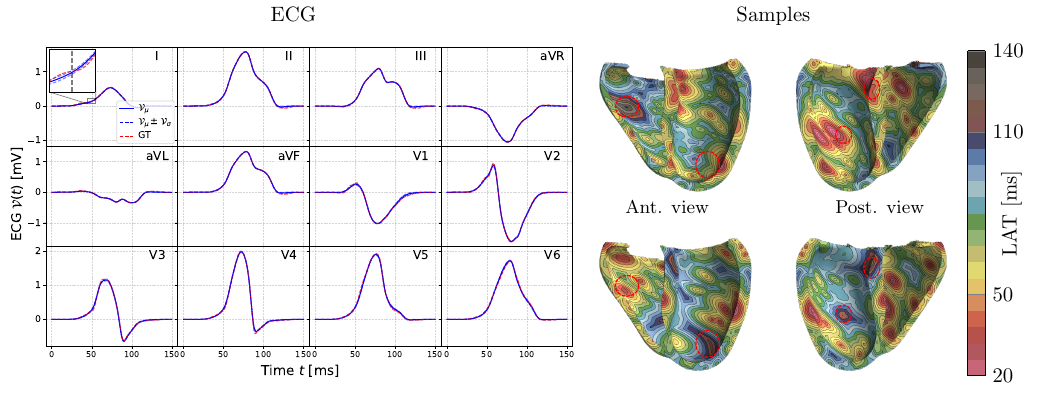}

    \caption{Identifiability problems of inverse \gls{ecg} methods. 
    Left panel: Using \geodesicBP{}, 
    the identification of optimal \glspl{pmj} locations and timings 
    to fit a given \gls{ecg} (\gls{gt}) can be achieved with high fidelity. 
    The \gls{ecg} distribution obtained from 20 optimization runs with different initializations
    forms a tight envelope around the ECG to be reconstructed (left, $\Vmean \pm \Vstd$).
    Right panel: The two samples with the largest difference in \gls{lat} \revaat{are shown
    that reveal significant variability} in the solution space which is indicative of limited identifiability
    in the general unrestricted case. 
    We highlight a few of the major differences in \gls{lat} (red dashed ellipses).
    This prompts for quantification of this variability and for \emph{a-priori} constraints 
    to reduce the non-uniqueness of the solution space.
    }
    \label{fig:ill_posed_heart}
\end{figure}
%


Optimization results using \geodesicBP{} for fitting the \gls{ecg} are shown in Fig.~\ref{fig:ill_posed_heart} (left panel) 
where the distribution \revaat{obtained} from $20$ optimized \glspl{ecg} (\Vmean{}, \Vstd{}) 
is overlaid on the \gls{gt} \gls{ecg} we aim to recover. 
The tight envelope around the measured \gls{ecg} demonstrates 
that \geodesicBP{} is able to fit the \gls{ecg} with very high fidelity \revaat{without scaling}. 
The maximum absolute error between optimized and \gls{gt} \gls{ecg} was $<\SI{2.8e-2}{\milli \volt}$ 
(relative error \SI[round-precision=2]{5.098930131663843}{\percent}) according to~\eqref{eq:ecg_dist} (average: \SI{2.07e-2}{\milli \volt}/\SI[round-precision=2]{4.066672407563584}{\percent}).
Similarly, the Pearson correlation coefficient of the \glspl{ecg} of all samples w.r.t.\ the \gls{gt} was $>\num{0.994}$.
Convergence is also rapidly achieved with \distecg{} falling below $\SI{1e-1}{\milli \volt}$ within $\le100$ iterations for most samples \revaat{(see also Fig.~\ref{fig:conv_plot})}.

However, the low errors observed in the \gls{ecg} did not necessarily correspond to small errors 
in the ventricular activation map $\tau$.
This is illustrated in Fig.~\ref{fig:ill_posed_heart} (right panel) 
where the two cases with the most extreme difference (\distlat{}) in the activation map are shown. 
The variability in reconstructed position of \glspl{pmj} was significant between these two samples, 
as readily evidenced by fundamentally different initial activation sites.
For instance, a site of earliest activation in the \gls{lv} mid-anterior endocardium in one run (top)
is a site of latest activation in another run (bottom).
Quantitatively, the absolute error in \gls{lat} was \SI{23.12}{\milli \second} 
(average between the samples: \SI{18.32}{\milli \second}).

\begin{figure}[tb]
    \centering
    \includegraphics[width=.9\linewidth, trim={0cm 4.65cm 0cm 0cm}, clip]{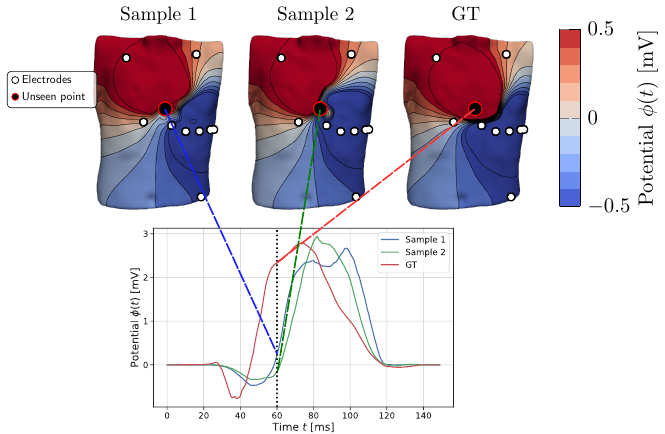}
    \caption{Shown are the two \glspl{bspm} of maximum distance among the $20$ samples 
    next to the sought-after \gls{gt} solution 
    at a single instant in time ($t=\SI{60}{\milli \second}$, top), 
    after normalization with respect to the mean unipolar potential (see Sec.~\ref{sec:pseudo_bidomain}). 
    Differences in potential fields are witnessed in space across the body surface, 
    as well as over time, when compared at a single unseen electrode location (bottom panel).
    At all \gls{ecg} electrode locations used for optimization no differences are apparent, 
    all potential traces are predicted with approximately zero loss (see Fig.~\ref{fig:ill_posed_heart}, left panel). 
    }
    \label{fig:ill_posed_torso}
\end{figure}

While the 12-lead \gls{ecg} was very similar across all samples, 
a more pronounced variability was observed in the \gls{bspm}, as shown in Fig.~\ref{fig:ill_posed_torso}. 
Again, we show the two \glspl{bspm} that were most different according to~\eqref{eq:pb_dist} (maximum absolute difference: \SI{0.11}{\milli \volt}), together with the \gls{gt} \gls{bspm}. 
It is worth noting that these two samples with maximum distance in \gls{bspm} shown in 
Fig.~\ref{fig:ill_posed_torso} are not the same samples 
that yielded a maximum distance in the activation maps shown in Fig.~\ref{fig:ill_posed_heart}.

Overall, the dipole-shaped torso potential was well-captured by all samples, 
indicating that the overall activation pattern was similar across all samples. 
The least correlated areas were found in the sternal region of the \gls{bspm}, 
where reconstructed \glspl{bspm} exhibited more complex patterns than the \gls{gt}. 
Selecting one point in this region shows different extracellular potentials between the two samples $\phi_i$, $\phi_j$ and the \gls{gt} $\phi_{\text{GT}}$ (bottom panel of Fig.~\ref{fig:ill_posed_torso}).
The correlation between the extracellular potentials at this unseen point was comparatively low,
with correlation coefficients against the \gls{gt} potentials of only $0.69$ and $0.64$, respectively.

\subsection{Identifiability under physiological constraints}
\label{sec:results_constr_bspm}

Unrestricted optimization with permissible \glspl{pmj} locations over the entire ventricular domain 
led to significant errors in the reconstruction of the ventricular activation sequence $\tau$, 
and to non-physiological activation sequences.
We investigated therefore whether imposing constraints based on physiological \emph{a priori} knowledge, 
as described in Sec.~\ref{sec:constraints}, can be effective in reducing the reconstruction error
and enhancing identifiability
\begin{figure}[htb]
    \centering
    \includegraphics[width=.93\linewidth]{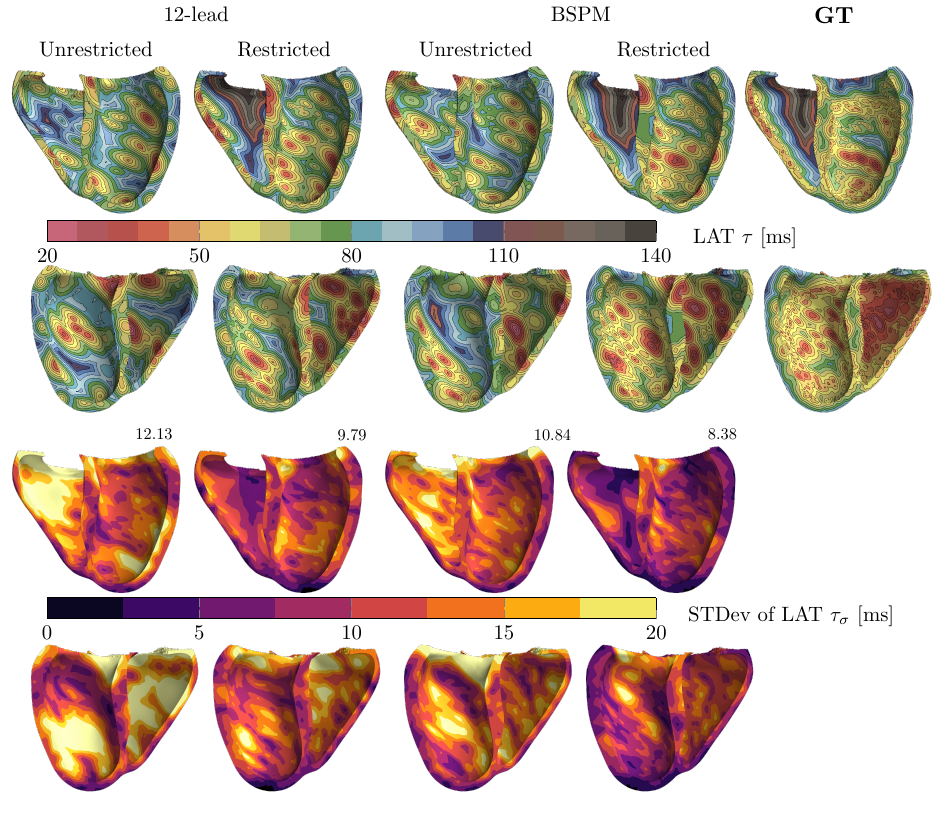}
    \caption{Preview of the distribution of the $20$ samples,
    reconstructed from 12-lead \gls{ecg} and \gls{bspm} with $128$ electrodes as observations, 
    both with and without physiological constraints.
    We visualize the sample most closely representing the mean $\taumean{}$ (top 2 rows) along with the deviation between the samples $\taustd{}$ (bottom 2 rows, see Sec.~\ref{sec:evaluation} for the exact definitions).
    The average deviation over the entire domain, $\bar{\tau}_\sigma$,
    is indicated at the top right corner of each $\taustd$ model.}
    \label{fig:results_mean_std}
\end{figure}

\revaat{In Fig.~\ref{fig:results_mean_std}, the model closest to the mean $\taumean$ (top) and variability models $\taustd$ (bottom) are shown.
We show the 12-lead and 128-electrode \gls{bspm} both in the restricted and unrestricted settings}. 
With restrictions, the average error in endocardial activation improved relative to the \gls{gt}, 
particularly in the \gls{rv} inferior wall. 
Non-physiological early activation in the basal and epicardial regions as seen in Fig.~\ref{fig:ill_posed_heart}  
is partially mitigated when restrictions are enforced, and more consistent with \gls{gt}, resulting in a more physiological transmural activation pattern
and a later activation of basal regions.


Imposing physiological restrictions and increasing the density of observations 
effectively reduce regions of pronounced variability.
This is reflected in an average deviation $\bar{\tau}_\sigma$, 
indicated in the top right corner of the variability models, 
which is always lower in the restricted setting.
Deviations ranged from \SI{12.13}{\milli \second} in the unrestricted case with only a 12-lead \gls{ecg} as observation 
down to \SI{8.38}{\milli \second} \revaat{in the restricted case} with a high density \gls{bspm}. 
Interestingly, an excellent match of the surface \gls{ecg} was reached even though
the physiological restrictions imposed during the reconstruction
did not exactly coincide with the prescriptions made during the generation of the \gls{gt} model. 
A geometric cutoff value of \SI{2.5}{\milli \meter} in the \gls{gt} model, as defined in Sec.~\ref{sec:constraints}, resulted in \SI{86}{\percent} capture of the \gls{pmj} nodes, with those exceeding \SI{2.5}{\milli \meter} all residing in the RV septum and free wall. A maximum penetration value of \SI{4.3}{\milli\meter} was reached, well above the cutoff value of \SI{2.5}{\milli \meter} used during the optimization from the \gls{ecg}.

\subsection{Identifiability with higher density observations}
\label{sec:dense_bspm}

\begin{figure}[htb]
    \centering
    \includegraphics[width=.8\linewidth]{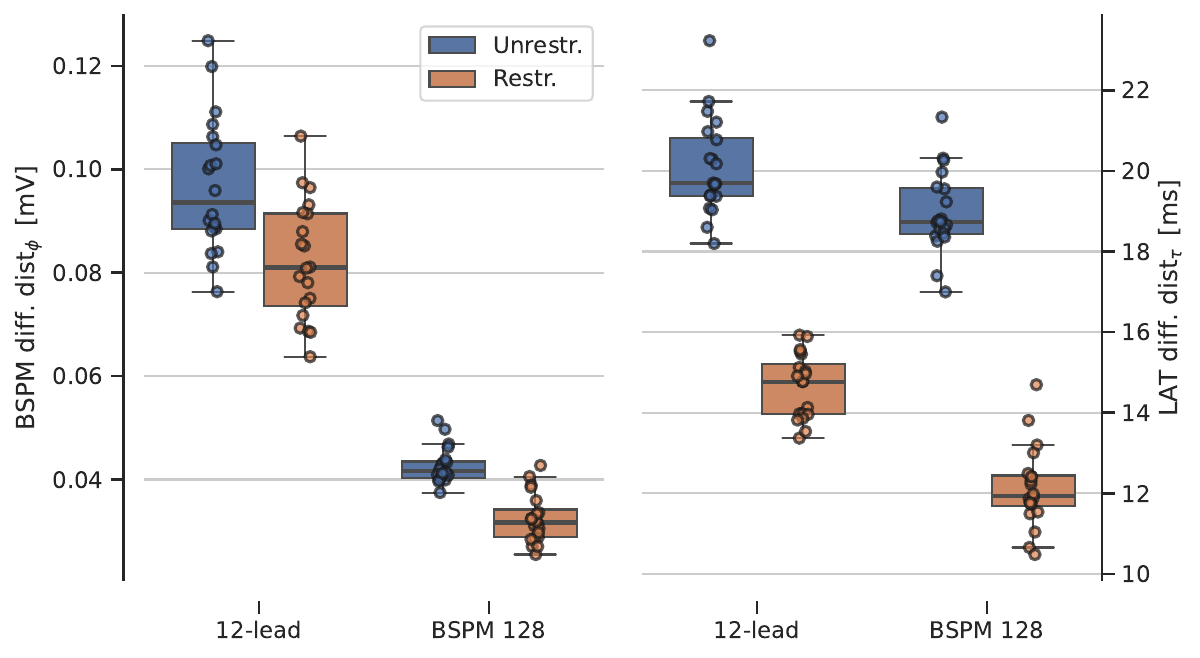}
    \caption{Comparison between pseudo-bidomain solutions on the torso 
    in terms of errors in potential reconstruction \distpb{}~\eqref{eq:pb_dist} (left panel) 
    and in the ventricular activation sequence \distlat{}~\eqref{eq:sol_dist} (right panel).
    Each pair of boxplots shows the errors for unrestricted (blue) and restricted (brown) cases, respectively.
    Higher density observations led to a significant improvement in the reconstruction of the \gls{bspm},
    but impact on the reconstruction of the ventricular activation map $\tau$ was marginal.}
    \label{fig:pb_comparison}
\end{figure}

Observation in Sec.~\ref{sec:uniqueness} indicates a perfect match at all \gls{ecg} electrode locations, 
and an overall close agreement in the reconstructed \glspl{bspm}, 
but deviations are witnessed for unseen positions on the torso 
where reconstructed and \gls{gt} \gls{bspm} do not match.
We investigate therefore whether an increase in electrode density on the torso is a suitable remedy 
to further improve \gls{bspm} reconstruction and identifiability of the \gls{lat} map $\tau$,
as suggested by the achieved reduction in regions of high variability in Fig.~\ref{fig:results_mean_std}.
%
Using higher density \gls{bspm} observations led to a considerable reduction in the reconstruction error 
over the torso, from \SI{9.67e-02}{\milli \volt} down to \SI{4.27e-02}{\milli \volt} on average (see Fig.~\ref{fig:pb_comparison}, left panel). 
However, this was less effective for improving the reconstruction of the ventricular activation sequence $\tau$
where only marginal improvements were witnessed, 
with errors decreasing from \SI{20.08}{\milli \second} down to \SI{18.95}{\milli \second} on average 
(Fig.~\ref{fig:pb_comparison}, right panel). 
As such, imposing physiological constraints is much more effective 
in reducing the reconstruction error in terms of \distlat{} than increasing the density of observations.
With physiological restrictions an average \distlat{} error of \SI{14.63}{\milli \second} and \SI{12.14}{\milli \second} was achieved for 12-lead \gls{ecg} and \gls{bspm} reconstructions, respectively.

\begin{figure}[htb]
    \centering
    \includegraphics[width=.8\linewidth]{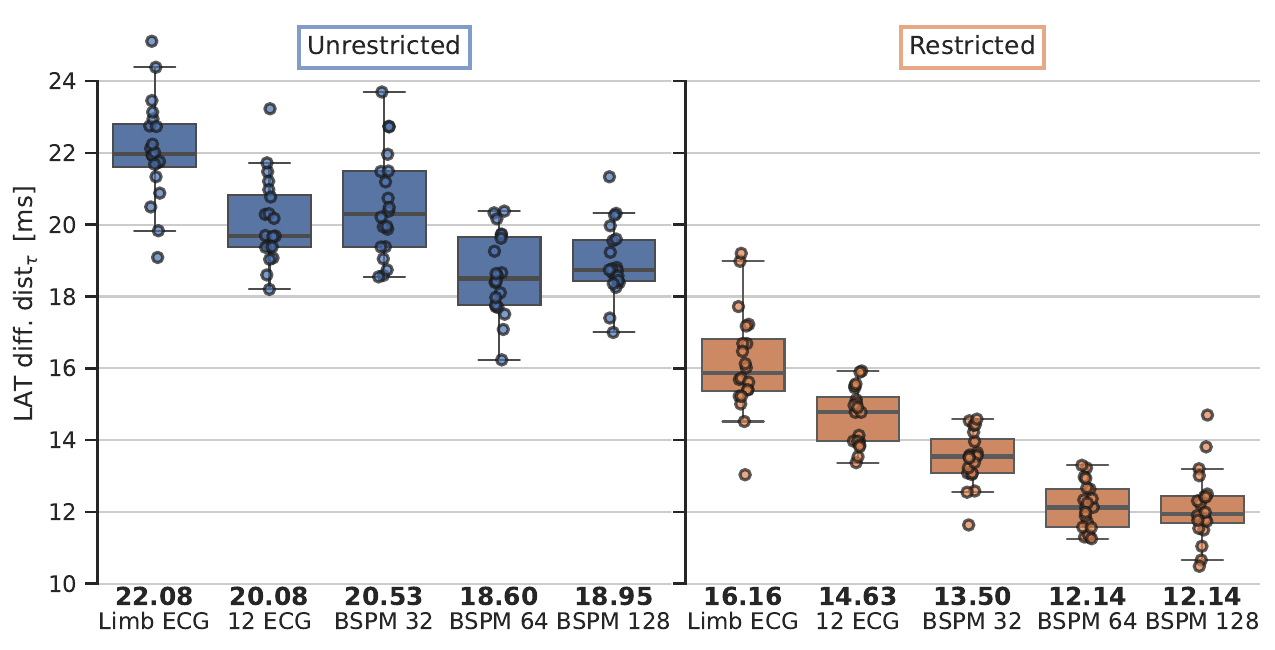}
    \caption{Comparison between \glspl{lat} as reconstructed with an increasing number of torso electrodes, both for ignoring (left panel) and enforcing (right panel) physiological restrictions.
    For each boxplot, the average \distlat{} with respect to the \gls{gt} of the 20 samples is shown as a number at the bottom line.}
    \label{fig:elec_comparison}
\end{figure}

A further differentiation of the impact due to observation density versus physiological constraints 
is shown in Fig.~\ref{fig:elec_comparison}.
The quality of the reconstructions of the ventricular activation sequence (\distlat{}) is compared 
for increasing observation density from 4 (limb leads) to 128 recording sites (\gls{bspm} vest),
comprising limb-lead and 12-lead \gls{ecg} as well as \gls{bspm} vests with increasing electrode density (Sec.~\ref{sec:setup_bspm}) with and without physiological constraints.
Increasing the number of electrodes improves the reconstruction of the activation map $\tau$, 
as the error in \distlat{} is decreased from an average \SI{22.08}{\milli \second} down to \SI{18.60}{\milli \second}.
In addition to the number of electrodes their placement is also a relevant factor,
as indicated by the slightly better reconstruction from the 12-lead \gls{ecg} 
over the higher density 32-electrode grid-based \gls{bspm} vest
where the placement of the 12-lead \gls{ecg} electrodes is denser in the proximity of the heart 
than the  strict grid-based 32-electrode \gls{bspm} vest.
Further increasing the density of electrodes of the \gls{bspm} vests past $64$ electrodes 
is only marginally beneficial for the reconstruction error, at the cost of slightly larger
variability in the samples.
However, as noted before, the major improvement in terms of \distlat{} results 
from imposing physiological constraints as seen in Fig.~\ref{fig:pb_comparison}.
With these the error in \distlat{} decreases with electrode density 
from \SI{16.16}{\milli \second} down to \SI{12.14}{\milli \second}.

%
%

\subsection{Hyperparameter sensitivity - the number of PMJs} 
\label{sec:nr_eas_study}
\begin{figure}[htb]
    \centering
    \includegraphics[width=.95\linewidth]{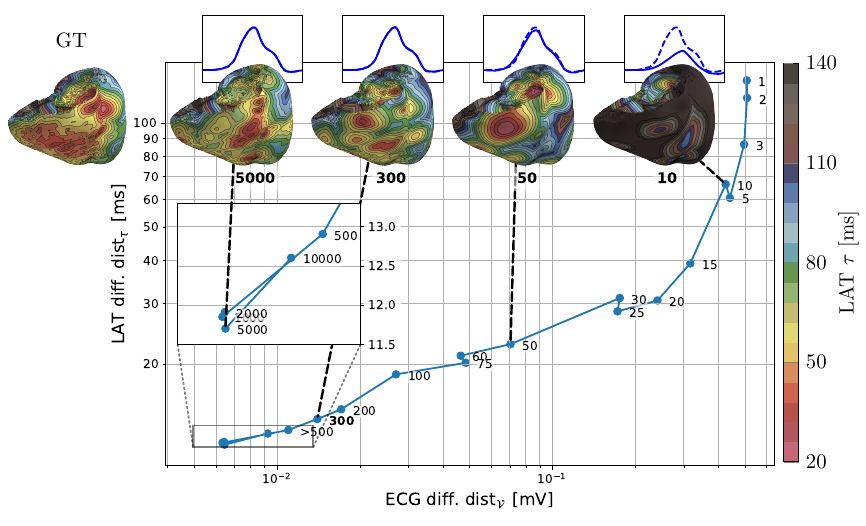}
    \caption{Hyperparameter study of the number of \glspl{pmj}.
        The 12-lead \gls{ecg} was fitted using a variable number of initial \glspl{pmj}
        and the error in reconstructing the \gls{ecg} and the activation map $\tau$ was measured against the \gls{gt} model
        in terms of \distecg{} and \distlat{}.
        A zoom-in on the results for $N \ge 500$ is shown in the inset.
        In the top panel optimization results are shown for a number of \glspl{pmj} of $N \in \{10, 50, 300, 5000\}$,
        comprising the reconstructed \gls{lat} maps 
        and the reconstructed Einthoven II \gls{ecg} (solid)
        next to the \gls{gt} \gls{ecg} (dashed).
        }
    \label{fig:nr_eas_study}
\end{figure}

Results obtained with optimization algorithms are sensitive to the choice of hyperparameters. 
In \geodesicBP{} 
-- since we are only using the \gls{ecg} in this study as a target in~\eqref{eq:inverse_ecg} -- 
the only hyperparameter to be chosen is the number of \glspl{pmj}.
These determine the number of earliest activation sites in the ventricles and, as such, 
the ability to approximate the ventricular activation map $\tau$.
The effect of the number of \glspl{pmj} upon the fidelity of reconstruction of 
the \gls{ecg} and the activation map $\tau$ was investigated in terms of \distecg{} and \distlat{}
by fitting to the \gls{ecg} with a variable number of \glspl{pmj}, 
ranging from 1 to 5000, with imposed physiological constraints.
Results are summarized in Fig.~\ref{fig:nr_eas_study}.

When using too few \glspl{pmj} ($<20$), 
the optimizer was unable to achieve a match of sufficient fidelity,
neither in terms of the \gls{ecg} (correlation $<0.8$, \distecg{} $>\SI{2.4e-1}{\milli \volt}$)
nor the ventricular activation (\distlat{} $>\SI{3e1}{\milli \second}$).
From 100 to 500 \glspl{pmj}, differences in terms of \distecg{} between the models essentially disappeared:
With a difference of $\SI{2.7006e-2}{\milli \volt}$ and a correlation of $0.998$ 
visually identical reconstructions of the 12-lead \gls{ecg} were obtained.
For even higher numbers of \glspl{pmj} $N > 500$, errors further improved in fitting \gls{ecg},
but only with marginal gains in \gls{lat} errors. Furthermore, the optimization
was highly susceptible to perturbations in the \gls{ecg} and the initial selection of the \glspl{pmj} prior
the optimization process. Overall, several \glspl{pmj} had no or low contribution to the activation.
With an initial $N = 300$ \glspl{pmj}, only \SI{23}{\percent} (or $\approx 70$) of them were responsible for activating \SI{95}{\percent} of the myocardium. For comparison, with $N = 5000$ initial sites, this number dropped to \SI{10.9}{\percent} ($\approx 545$ ) at the end of the optimization process. Inactive \glspl{pmj} for $N=300$ and $N=5000$ were \SI{62.33}{\percent} and \SI{78.5}{\percent}, respectively.
These numbers were calculated using the \gls{roi} defined in Sec.~\ref{app:roi}.

\section{Discussion}
\label{sec:discussion}

In this work, we have demonstrated the ability of a novel approach 
for constructing \glspl{cdt} of human ventricular electrophysiology 
with unprecedented fidelity and speed using only non-invasive clinical data 
-- the 12-lead \gls{ecg} and tomographic images. 
The superiority of our approach is based on its efficiency, robustness, and accuracy 
in reproducing clinically recorded signals \revaat{without scaling} that cannot be matched by \revaat{any other} current state-of-the-art \revaat{methodology. 
Our} findings indicate that incorporating histo-anatomical knowledge 
on the ventricular conduction system into our fitting approach 
facilitates an accurate recovery of \gls{lat} only with the clinical 12-lead \gls{ecg}\revaat{.
Even in the most complex scenario of the Purkinje-mediated ventricular activation during normal sinus rhythm, there is limited need for high-density \gls{bspm} data that are, in general, clinically not available.} 
This novel optimization technique represents a significant advancement in generating credible \glspl{cdt} \revaat{ 
in a scalable manner and with demonstrable fidelity. 
Our approach holds considerable promise for translating \glspl{cdt} into clinical practice.}


We provide evidence that \glspl{cdt} cannot be uniquely identified from surface
\gls{ecg} recordings, even when using high-density \gls{bspm}.
\revaat{In fact}, multiple markedly distinct activation patterns can
predict the same surface \gls{ecg}, within a tolerance level comparable to
noise.
As a consequence, quantifying the uncertainty in the reconstruction process is crucial for clinical applications. 
In this study, we show that the variability in reconstruction can be significantly reduced by imposing physiological constraints on the permissible domain for \glspl{pmj}. 

Our study further proved that variability in the \revaat{\gls{lat}} reconstruction is more effectively reduced 
by using \emph{a priori} knowledge \revaat{in the form of} constraints than by increasing the density of observations.
Specifically, when constructing a \gls{cdt}, high-density \gls{bspm} offers only marginal benefits 
compared to the standard 12-lead \gls{ecg}. 
Although this may seem counter intuitive at first, our analysis revealed 
that in a population of \glspl{cdt} all \revaat{fitted} to the same 12-lead \gls{ecg}, 
the residual variability in the \gls{bspm} is predominantly confined to a narrow region around the sternum. 
Accurately resolving the spatio-temporal evolution of potential fields in closer proximity of the heart,
including non-dipole components of the field, would require 
advanced \gls{bspm} systems capable of higher spatial resolutions. 


Our method recovers the ventricular conduction system by optimizing the location and timing of \glspl{pmj} within a subendocardial domain -- 
defined on the basis of histo-anatomical \emph{a priori} knowledge. 
In contrast, other methods focus on modeling the entire \gls{hps}, 
which can be personalized to match 12-lead \gls{ecg} or \gls{eam} data~\cite{camps2024digital,alvarez_barrientos_probabilistic_2025,barber2021estimation}. 
Given the histological evidence and current understanding of the \gls{hps} (see Sec.~\ref{sec:gt_model}), 
it is reasonable to generate it using a procedural approach, 
such as tree generation algorithms~\cite{costabal2016generating} or constrained optimization methods~\cite{berg2023enhanced}. 
In these approaches, the position and timing of \glspl{pmj} are fully determined by the tree's topology and the assumed \gls{cv}. 
No such Purkinje restriction is in place for the \glspl{pmj} created with our method,
thus allowing for a more flexible placement and simpler manipulation of \glspl{pmj} 
as compared to a topological representation of a Purkinje network.
However, for a fully mechanistic representation in a \gls{cdt} an explicit topological model
of the \gls{hps} is required 
to facilitate an accurate prediction of ventricular activation and \glspl{ecg} under other conditions 
the \gls{cdt} was not calibrated for such as pacing.
In our approach a suitable Purkinje network that produces an equivalent ventricular activation sequence  
must be constructed for a given set of \glspl{pmj} in a post-processing step after optimization.
The construction of such \gls{hps} models from \glspl{pmj} is readily achieved for lower dimensional approximations 
of the \gls{hps} such as assuming a five fascicular \gls{hps} to activate a fast conducting subendocardial layer,
as has been shown previously \cite{gillette2021hps}, 
but a robust reconstruction of a \gls{hps} to match the fine-grained activation pattern as produced with our approach, 
with a number of \glspl{pmj} $N \ge 5$ has not been reported yet.

A crucial differentiating factor of our method is the use of gradient information
in the fitting procedure.
Other methods lack this information which necessitates computationally expensive sampling techniques to optimize the QRS complex of the \gls{ecg}. 
Even when using the fastest \gls{ecg} forward models 
able to predict \glspl{ecg} with close to real-time performance~\cite{neic2017efficient} 
or even faster~\cite{pezzuto_evaluation_2017},
these methods often require many hours to compute or result in sub-optimal \gls{ecg} fits. 
Surrogate models based on Gaussian process regression~\cite{pezzuto_learning_2022} or neural networks~\cite{salvador_digital_2023} are promising options to further increase by one or two orders of magnitude the forward solver, but their construction is inherently tied to sampling in the parameter space. Furthermore, defining surrogate models in generic geometries is still an open challenge.

\begin{figure}[t]
    \centering
    \includegraphics[width=.85\textwidth]{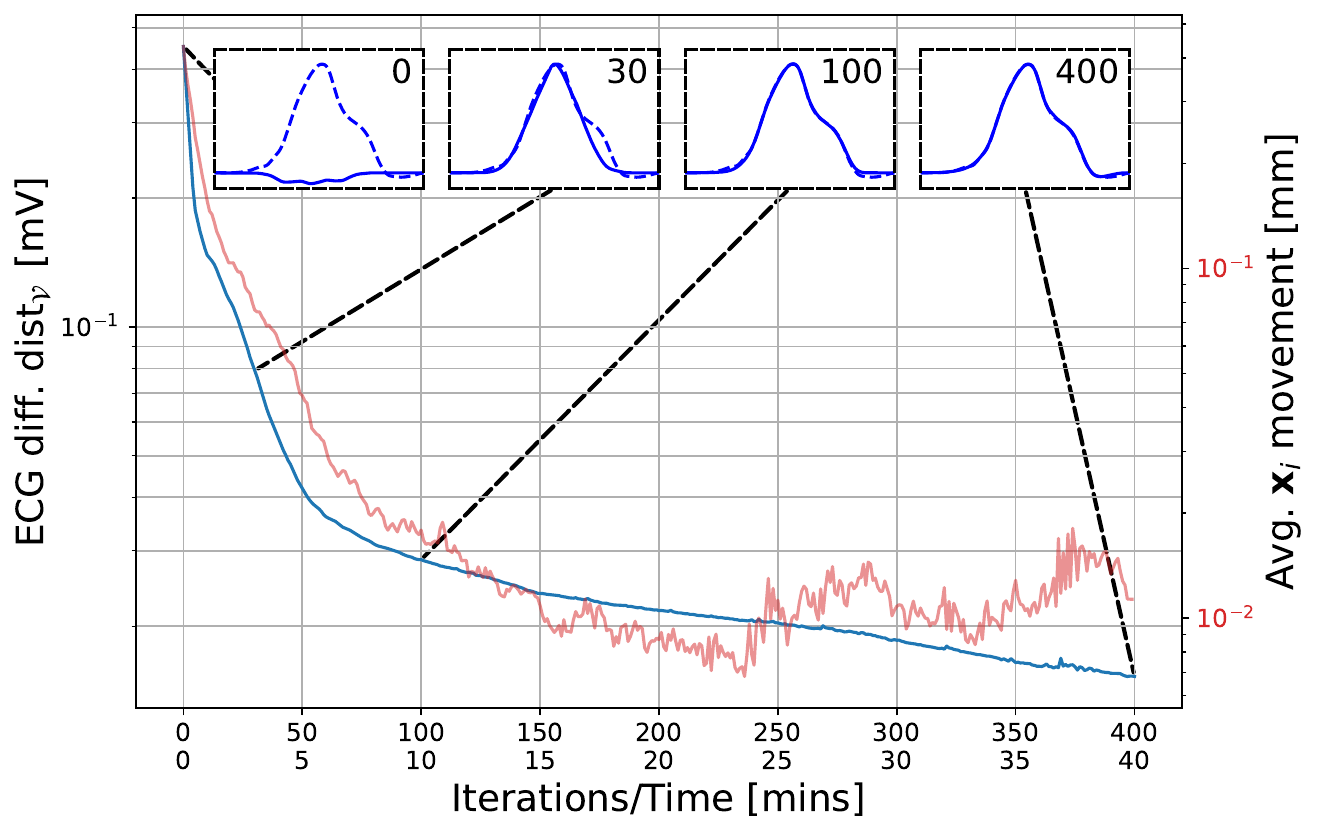}
    \caption{Convergence history of Geodesic-BP:
    We show the evolution of the solution in terms of the \gls{ecg} on the left y-axis (logarithmic), over the iterations (first row) and approximate computational time (second row) on the x-axis.
    The top $4$ plots show the computed lead II \gls{ecg} in the current iteration (solid)
    next to the \gls{gt} \gls{ecg} (dashed).
    Additionally, we show the averaged movement of the initial sites $\vx_i$ w.r.t.\ the last iteration in red (right $y$-axis, also in logarithmic scale).}
    \label{fig:conv_plot}
\end{figure}

The performance of technologies for calibrating \glspl{cdt} is a key factor 
determining the application scope.
Our method is computationally expensive, 
but the implementation is heavily optimized to take advantage of \gls{gpu} computations,
thus making it very fast in comparison to other similar methods.
The absolute execution time of each run depends mostly on the longest geodesic path in the mesh, which is related to the number of active \glspl{pmj} and the time needed to compute the projections.
All samples for $N=300$ take less than \SI{30}{\minute} on a NVidia RTX 4090 
for a single optimization run of 400 iterations as described in Sec.~\ref{sec:geodesic_bp}. 
\revaat{As shown in Fig.~\ref{fig:conv_plot}}, 100 iterations are often sufficient to achieve a good fit of the \gls{ecg}.
\revaat{To note, the ``elbow'' point in the \gls{ecg} fit is around $100$ iterations, and also the movement of the \glspl{pmj} decrease rapidly, partly because many sites deactivate, but also since the gradient magnitude decreases.}
In general, a single iteration step takes less than \SI{6}{\second},
which is comparable to the cost of running a forward simulation.

%
%


For all experiments conducted in this study the initial number of \glspl{pmj} was chosen as $N=300$, 
justified by the insights gleaned from the hyperparameter study in Sec.~\ref{sec:nr_eas_study}. 
Since the \glspl{pmj} are \revaat{uniformly} initialized over the \revaat{entire time range and admissible domain}, 
a large fraction of \glspl{pmj} is inactive at the beginning of the optimization.
The number of \glspl{pmj} suggested in Sec.~\ref{sec:nr_eas_study} may appear excessively large 
in comparison to other studies~\cite{gillette2021framework, grandits_inverse_2020,cardone2016human,camps2024digital}
that report only 5-7 \glspl{pmj}.
However, this is readily explained by a fundamental difference in the underlying model assumption
that builds on a fast-conducting subendocardial layer, aiming to mimic the fast isotropic activation 
mediated by the terminal Purkinje network of the \gls{hps} \cite{myerburg1972physiology}. 
Similar studies using geodesics report even higher numbers \revaat{of \glspl{pmj}} to achieve a good match 
with endocardial electro-anatomical mapping data \cite{cardenes2015estimation}, 
and studies that do not rely on fast-conducting layers report numbers in the range of 11-51 \glspl{pmj} \cite{barber2021estimation}. 
This is in agreement with our findings,
showing the requirement of at least $\ge 50$ \glspl{pmj} to effectively match a given QRS complex in the \gls{ecg} in the absence of a fast conducting layer. 



The most important factor for achieving a high fidelity reconstruction of both 
\gls{ecg} and activation sequence was the usage of physiological constraints 
detailed in Sec.~\ref{sec:constraints} that greatly improved the reconstruction quality. 
While exact knowledge on the structure of the \gls{hps} could be, in principle, prescribed 
according to the \gls{gt} model, this is not the case in real world applications 
when \glspl{cdt} are built from patient data. Applying physiological constraints to real patient data is not straightforward 
as the exact structure of the \gls{hps}, 
as well as the distribution and behaviour of \glspl{pmj} within the ventricular, 
is not well understood in the wider population\revaat{. Even less is known on the inter-individual variability}. To elucidate the effects that unknown or improper assumptions of the \gls{hps} may play in the ability to obtain an exact match in the \gls{ecg}, exact constraints from the \gls{gt} model were not utilized. Namely, the coverage of endocardial layer $S_e$ during optimization coincided with the \gls{gt} model, but a singular \gls{pmj} penetration depth $d_\text{pmj}$ of \SI{2.5}{\milli \meter} was deployed.  Regardless of the generalized assumption in the \gls{pmj} penetration depth $d_\text{pmj}$, 
a perfect match could be obtained for the \gls{ecg} of the simulated \gls{gt} model exhibiting robustness in the optimization method. 

Further exploration of alternative geometrical constraints in real world applications when \glspl{cdt} are built from patient data is thus needed. 
Alternative values for \gls{pmj} depth could be prescribed as ventricular wall thicknesses are highly variable within the healthy population, but range between 4-\SI{12}{\milli \meter} in the LV \cite{sjogren1971left} and 2-\SI{5}{\milli \meter} in the \gls{rv} \cite{prakash1978determination} during diastole. 
Furthermore, sparsity of the \gls{hps} in the RV inferior wall was deployed based on histological studies \cite{stephenson_2017_high}, but this may not agree with full coverage of the RV inferior wall with Purkinje fibers \cite{myerburg1972physiology}.
\revaat{A possibility is also to reformulate the inverse problem in a Bayesian framework, with a strong prior distribution for the \gls{hps}, for instance based on rule-based or generative methods~\cite{alvarez_barrientos_probabilistic_2025}.}

Moreover, at the discrete level of the advanced \gls{gt} model used in this study\revaat{,}
the information on the \gls{pmj} depth is spatially highly variable\revaat{. 
The terminal} \glspl{pmj} of the Purkinje network are resistively coupled to the ventricular myocardium 
within a prescribed radius and minimum number of connections\revaat{. This accounts} for electrotonic loading and associated anterograde and retrograde conduction delays \cite{boyle10:_purkinje_defib}. 
With the coarser mesh resolution of around \SI{1.2}{\milli \meter}\revaat{, a spatial }
variability at the scale of the mesh resolution is inevitable.

This paper only focused on the optimization of \glspl{pmj} for 
a varying number of electrodes, but always assumed knowledge of the anatomy, lead fields and \gls{cv} tensor.
However, even in this restricted scenario, optimization of only the \glspl{pmj} remains an ill-posed task 
as we have shown that requires strong constraints even for the already physics-restricted eikonal model.
A further extension to optimize the \gls{cv} tensor is possible from a methodological point-of-view~\cite{grandits_inverse_2020}, 
but is guaranteed to further expand the space of possible solutions and would require additional regularization considerations. For instance, additional information may be based on imaging or \gls{eam} data, which can be used to reconstruct the anisotropic structure of the tissue~\cite{Grandits2021PIEMAP,Lubrecht2021PIEMAP, RuizPINN2022}.
A further limitation \revaat{is} the various anatomical simplifications inherent in the process of anatomical model generation
where structures pertaining to the non-compact myocardium 
including trabeculation, papillary muscles\revaat{,} or the moderator band are assigned to the blood pool.
As such, a significant discrepancy arises between the relatively smooth endocardial surfaces in the model
and the highly structured and rugged surfaces \emph{in vivo} in a patient.
Moreover, these ignored structures are tightly linked to the function of the \gls{hps}. 
For instance, the root regions of papillary muscles are activated by the Purkinje network, 
and the moderator band in the \gls{rv} encloses a fascicle of the \gls{hps}~\cite{gettes_effects_1968,walton_compartmentalized_2018}.
These are major confounding factors as the optimization of \gls{pmj} locations is carried out 
within a sub-endocardial space that may be strikingly different between model and \emph{in vivo} in a patient.
Building such highly detailed anatomical models that accurately represent the non-compact myocardium is feasible~\cite{plank09:_generation} and can be achieved also with routine clinical imaging.
These effects require further investigation, but this was beyond the scope of this study.

\revaat{
In this work, we have not considered the repolarization phase and the corresponding T-wave. In principle, our methodology can be easily extended to optimize specific repolarization parameters in the forward model~\cite{gillette2021framework, camps_harnessing_2025}.
}

In conclusion, in this paper we extensively validated our recent method \geodesicBP{} 
for creating \gls{cdt} from the surface \gls{ecg} 
and discussed the identifiability and uniqueness in recovering \gls{lat} maps.
While the \gls{ecg} reconstruction can be achieved in a reliable manner, 
the quality of the \gls{lat} reconstruction highly depends on the restrictions put onto the \glspl{pmj}, the initialization of the parameters, and the number of electrodes. 
According to our findings, the restrictions were the most important factor to improve the overall quality of the results and minimize variability between the reconstructions.










\section*{Acknowledgements}
SP and GP acknowledges the support of the SNF-FWF ``CardioTwin'' project (grant number 214817). SP also acknowledges the PRIN-PNRR project no.~P2022N5ZNP and INdAM-GNCS.
The computational results presented have been achieved in part using the \gls{vsc} (project no.~71138, TG) and the CSCS-Swiss National Supercomputing Centre project no.~s1275 (SP).

\section*{Declarations}

Author GP holds shares in NumeriCor \emph{GmbH} but declares no non-financial competing interests.
All other authors similarly declare no non-financial competing interests.

\ifpreprint
\section*{Code and Data availability}

The datasets generated and analyzed for the 12-lead case (restricted and unrestriced) will be made available upon acceptance of the manuscript.
Similarly, the code to compute the \gls{lat} from the \glspl{pmj} using \texttt{fim-python}~\cite{grandits_fast_2021} and the code for the subsequent computation of the 12-lead \gls{ecg} for the restricted and unrestricted case will be provided upon acceptance.
A confidential, temporary preview version of the code and dataset is available to the reviewers online\footnote{\url{https://cloud.uni-graz.at/s/gBW3CHYkFB5TNMM}}.
\fi

\section*{Author contributions}

\textbf{Thomas Grandits}: Conceptualization, Methodology, Software, Validation, Formal analysis, Investigation, Resources, Data curation, Writing -- original draft, Writing -- review \& editing, Visualization.\\
\textbf{Karli Gillette}: Conceptualization, Funding acquisition, Formal analysis, Investigation, Methodology, Validation, Writing -- original draft, Writing -- review \& editing, Project administration.\\
\textbf{Gernot Plank}: Conceptualization, Resources, Writing -- original draft, Writing -- review \& editing, Project administration, Funding acquisition.\\
\textbf{Simone Pezzuto}: Conceptualization, Resources, Writing -- original draft, Writing -- review \& editing, Project administration, Funding acquisition.


\bibliographystyle{elsarticle-num-names}
\bibliography{biblio}

\appendix

\section{Region of influence (ROI)} 
\label{app:roi}
In order to optimize problem in~\eqref{eq:inverse_ecg} we compute in each iteration of the optimization algorithm the solution $\tau$ and its gradient $\nabla _{(\vx_i, t_i)} \tau$ for each \gls{pmj} $(\vx_i, t_i)$ in one backpropagation.
The gradient is used in the optimization to match the \gls{ecg}, but doubles as a way to compute the influence of a \gls{pmj} on the total solution.
To better understand the concept, consider the \gls{lat} at a point $\vx \in \Omega$, which could be equally defined as the distance \revaat{(in geodesic terms) from the closest initial site} to the point $\vx$, i.e.
\begin{equation}
    \tau(\vx) \coloneqq \min_i t_i + \int_0^1 \gamma_i(t) \diff{t}
    \label{eq:characteristic}
\end{equation}
where $\gamma_i$ is the geodesic path between $\vx_i$ and $\vx$.
Using~\eqref{eq:characteristic}, the derivative of each point w.r.t. $t_i$ readily follows:
\begin{equation}
    \frac{\partial \tau}{\partial t_i} (\vx) =
    \begin{cases}
        1 & \text{if } (\vx_i, t_i) \text{ activates } \vx \\
        0 & \text{else}
    \end{cases}
    \label{eq:characteristic_diff}
\end{equation}
In \gls{ep} terms,~\eqref{eq:characteristic_diff} can be interpreted as finding for each point $\vx$ in the domain its closest initial site $(\vx_i, t_i)$ (excluding equidistant boundary regions with measure $0$), which is responsible for the activation of $\vx$.
Using~\eqref{eq:characteristic_diff}, we can give a mathematically and logically sound quantitative definition of the \gls{roi} for each initial site $(\vx_i, t_i)$:
\begin{equation}
    \text{ROI}_i = \int_\Omega \frac{\partial \tau}{\partial t_i} (\vx) \diff{\vx},
    \label{eq:roi}
\end{equation}
which is measuring the volume of tissue that is activated by the initial site $(\vx_i, t_i)$.
The computation of all \glspl{roi} can be performed using a single backpropagation of $\int_\Omega \tau(\vx) \diff{\vx}$.
Note that initial sites with $\text{ROI}_i = 0$ can exist, which just means that they do not activate any tissue and will not be optimized by the gradient-based optimization.
This is also consistent with the fact that their presence or absence does not change the computed solution $\tau$.
We call an initial site $(\vx_i, t_i)$ with $\text{ROI}_i = 0$ simply an \emph{inactivate initial site}.
Further information on the relation between the method of characteristics and the eikonal equation can be found in~\cite{carmo_riemannian_1992,grandits_geasi_2021}.

\end{document}